\documentclass[11pt]{article}
\usepackage{amssymb,latexsym,amsfonts,verbatim,amscd}

\newtheorem{Def}{Definition}[section]
\newtheorem{Thm}[Def]{Theorem}
\newtheorem{Lem}[Def]{Lemma}
\newtheorem{Prop}[Def]{Proposition}
\newtheorem{Cor}[Def]{Corollary}
\newtheorem{Rem}[Def]{Remark}

\def\telos{\hfill$\dashv$}

\begin{document}
\sloppy

\title{Localizing the axioms}
\author{Athanassios Tzouvaras}

\date{}
\maketitle

\begin{center}
 Department  of Mathematics\\  Aristotle University of Thessaloniki \\
 541 24 Thessaloniki, Greece \\
  e-mail: \verb"tzouvara@math.auth.gr"
\end{center}

\begin{abstract}
We examine  what happens if we replace ZFC with a
localistic/relativistic system, LZFC,  whose central new axiom,
denoted by $Loc({\rm ZFC})$,  says that every set belongs to a
transitive model of ZFC.  LZFC consists of $Loc({\rm ZFC})$ plus
some elementary axioms  forming Basic Set Theory (BST).  Some
theoretical reasons for this shift of view are given. All $\Pi_2$
consequences of ZFC are provable in ${\rm LZFC}$.   LZFC
strongly extends  Kripke-Platek (KP) set theory minus $\Delta_0$-Collection
and minus $\in$-induction scheme.  ZFC+``there is an inaccessible cardinal'' proves the consistency
of LZFC. In LZFC we focus on models rather than cardinals, a
transitive model being considered as the analogue of an
inaccessible cardinal. Pushing this analogy further we define
$\alpha$-Mahlo models and $\Pi_1^1$-indescribable models, the
latter being the analogues of weakly compact cardinals. Also
localization axioms of the form $Loc({\rm ZFC}+\phi)$ are
considered and their global consequences are examined. Finally we
introduce the concept of standard compact cardinal (in ZFC) and
some standard compactness results are proved.

\end{abstract}


\section{Introduction}

The purpose of this paper is to look at ZFC from a certain
localistic/relativistic point of view. In current set theory we
believe that there is an objective reality of sets, the ``real
world'' $V$, the main properties of which are captured by  the
axioms of ZFC. In other words,  the ZFC axioms are supposed to
hold  {\em in} $V$. This is the absolutistic point of view. An
opposite  view, that may be called localistic/relativistic, would
consist in  claiming  that the ZFC axioms, especially  the
problematic  axiom of Powerset (and perhaps Replacement), should
refer not to $V$ itself but only  to  several {\em local} models,
which are counterparts of the reference frames of physics.
Conceivably there are more than one ways to formalize  this
general idea of local truth and  local models.  The formal account
presented in this paper is just one among them. Its main points
are roughly the following: (1) All local models of ZFC (or
extensions of it) that we consider are standard transitive sets.
(2) There is an abundance of them across the universe.

The motivation for such a shift of view comes from  the well-known
relativity, first pointed out by Skolem \cite{S22}, that occurs in
all first-order axiomatizations of set theory. Some fundamental
notions, especially cardinality and powerset, raise such
unsurmountable difficulties when treated as absolute entities,
that, until one comes up with a revolutionary new idea about what
the  powerset of an infinite set actually contains - which  possibly (though not necessarily) might
settle also the problem of  counting its members - one
would better let aside  the idea that ${\cal P}(\omega)$ exists in
$V$ and  instead be content with the idea that  ${\cal
P}(\omega)$ is a set with respect to transitive set-universes
{\em only,} i.e., in the local/relative form ${\cal P}^M(\omega)={\cal
P}(\omega)\cap M$, where $(M,\in)$ is a transitive  model of
ZFC.\footnote{Throughout the term ``transitive model'' is used
instead of the more cumbersome ``standard transitive (set)
model'', i.e., a transitive set $x$ equipped with the standard
membership relation $\in$, so that $(x,\in)\models{\rm ZFC}$. We
could just say  ``standard model'', would transitivity not be
independent from  standardness. A transitive set on the other hand
is implicitly thought as being structured by $\in$. In view of the
Mostowski's isomorphism theorem however, a standard model of ZFC
is essentially identical to a standard transitive one.} ${\cal
P}(\omega)$ itself makes sense only as a proper class. In
compensation  one may assume that transitive models of ZFC exist
{\em everywhere in $V$}, specifically that every set $x$ belongs
to some transitive model $y$. Such a view on the
one hand does not have any negative impact  on the study of
various kinds of infinite cardinals. For example it by no means
invalidates the theory of large cardinals, except of course that
these are now treated as relativized entities living only in
models. And on the other hand it spurs the interest in transitive
models themselves, as   objects of study {\sl per se} rather than  just
 a means. Large cardinals in particular constitute  a source of  ideas and
techniques some of which can be transferred  to models in order to
build analogous classifications among them.

A theoretical justification of the  above viewpoint is summarized
in the  following argument: Although we may believe that $V$ is
indeed an objective, absolute reality, it does not necessarily
follow that all properties and facts concerning objects of $V$
should be absolute too. {\em Some} properties may be subject
always and by their nature to local constraints, so that any
absolutistic judgment about them would simply not make sense. A
helpful and convincing analogy comes from the universe of physical
objects. According to the established paradigm of Relativity
Theory, this universe is also an objective, absolute reality of
things,\footnote{Or, at least, it {\em can} be. Obviously no final
decision  can be reached on such a metaphysical issue.} but
fundamental physical magnitudes like mass, length, time, velocity,
etc, are inherently relative quantities, heavily depending on the
observer's reference frame. If  fundamental attributes of physical
objects such as mass and size   are relative, why should the {\em
type } (or {\em degree}) {\em  of infinity} of an infinite set be
absolute? Of course there are differences: In the case of physical
universe there are experiments and measurements supporting  the
view of Relativity Theory, while for the universe of abstract sets
one can only make assumptions. Also one tends to accept  much more
easily that almost all physical properties (color, shape, smell,
etc) are subject  to relativization, than that this is also the
case with abstract properties, like number and structure, which
are commonly supposed to reflect deeper and more permanent
characteristics of beings. And in fact, finite cardinalities
$0,1,2,\ldots$  do not seem to relativize in any reasonable way.
But the various {\em infinite cardinalities} is a different
matter. Among all mathematical objects these should be the most
naturally expected to be inherently relative. A strong indication
is the ease by which the cardinality of  an infinite set can
change by means of forcing constructions.

So much for the viability of the localistic/relativistic
approach to set theory.  The purpose of the paper is to set out a
particular {\em implementation} of this approach through an
axiomatic system and examine  its logical strength and its set
theoretic consequences. The paper is organized as follows:

In section 2 we define the system LZFC (from ``local ZFC'') whose
main axiom is:
$$(Loc({\rm ZFC}))
\hspace{.5\columnwidth minus .5\columnwidth} \forall x\exists
y(x\in y \ \wedge \ Tr(y) \ \wedge \ (y,\in)\models{\rm
ZFC}).\hspace{.5\columnwidth minus .5\columnwidth} \llap{}$$ The
other axioms, forming the system BST (of Basic Set Theory), are
elementary assumptions like Pair, Union, etc, needed only  to
formulate $Loc({\rm ZFC})$. LZFC  proves all $\Pi_2$ consequences
of ZFC. Also LZFC proves the equivalence of the ${\rm Found}^*$ of $\in$-induction and the scheme ${\rm Found}_{On}$ of induction over the ordinals. However none of them seems to be derivable in LZFC. Consequently, transfinite induction along $On$ is not available in LZFC.
LZFC  does not prove  $\Pi_2$-Reflection, since  ${\rm
LZFC+\Pi_2}$-{\rm Reflection} \ $\vdash Con({\rm LZFC})$.  $\Sigma_1$-Collection is equivalent to $\Delta_0$-Collection over LZFC, but it is open whether the latter proves $\Delta_0$-Collection. The
class $L$ of constructible sets is definable (though one cannot
prove in LZFC that $L$ is an inner model of LZFC). Also standard
facts and constructions, like Completeness theorem,
L\"{o}wenheim-Skolem theorem, generic extensions, Mostowski
collapse etc, are available in LZFC. LZFC is a strong  extension of KP (Kripke-Platek set theory) minus
$\Delta_0$-Collection and minus the scheme ${\rm Found}^*$ of $\in$-induction. Concerning consistency, LZFC is a subtheory of ZFC+ ``there is
a proper class of inaccessible cardinals''. Also ${\rm ZFC}$ + ``there is an inaccessible cardinal'' proves the consistency
of ZFC+LZFC, while  ${\rm ZFC}$ + ``there is a natural model of ZFC'' proves the consistency of LZFC.

In section 3 we  discuss infinite (uncountable) cardinals and powersets (of infinite sets) in LZFC. In view of the absence of transfinite induction, no general statement about cardinals $\omega_\alpha$ and powersets  ${\cal P}^\alpha(\omega)$ can be derived.  Yet certain implications concerning    existence and absoluteness of concrete classes like  $\omega_1$, ${\cal P}(\omega)$ and $H(\omega_1)$ (and more generally  $\omega_n$, ${\cal P}^n(\omega)$ and $H(\omega_n)$, for $n\in \omega$) can be established. For example it is proved that   $H(\omega_1)\in M$ implies ${\cal P}(\omega)^M={\cal P}(\omega)$ and $\omega_1^M=\omega_1$; ${\cal P}(\omega)\in M$ implies $\omega_1^M=\omega_1$ and $H(\omega_1)^M=H(\omega)$, etc. We discuss also an ambiguity concerning the meaning of the symbols $\omega_\alpha$, for $\alpha>0$, and how it can be raised.

In section 4 we define  $\alpha$-Mahlo models  as analogues of
$\alpha$-Mahlo cardinals. This is a pretty natural notion: A model
$M$ of ZFC is Mahlo if the set of models of ZFC that belong to $M$
is a  stationary subset of $M$. Stationary, as well as closed
unbounded subsets of $M$, are restricted to {\em definable}
subsets of $M$.  Definability guarantees that the property of
$\alpha$-Mahloness is absolute for transitive models of ZFC. It is
shown in ZFC that if $\kappa$ is $\alpha$-Mahlo, then $V_\kappa$
is an $\alpha$-Mahlo model.

In section 5 we define in LZFC $\Pi_1^1$-indescribable models, as
analogues of $\Pi_1^1$-indescribable (i.e., weakly compact)
cardinals.  Concerning the existence of such models, we show (in
ZFC)  that if $\kappa$ is weakly compact then $V_\kappa$ is
$\Pi_1^1$-indescribable. Moreover, if $M$ is
$\Pi_1^1$-indescribable, then it is $\alpha$-Mahlo for every
$\alpha\in M$.

In section 6 we consider localization axioms of extensions of ZFC,
i.e., of  the form $Loc({\rm ZFC}+\phi)$, or $\{Loc({\rm
ZFC}+\phi):\phi\in\Gamma\}$, for some set of sentences $\Gamma$,
and examine their consistency (when added to LZFC)  and their impact on $V$. For instance it is shown that  for every $\Pi_1$ or $\Sigma_1$ sentence $\phi$,
$Loc({\rm ZFC}+\phi)+Loc({\rm ZFC}+\neg\phi)$ is
inconsistent. Further, $Loc({\rm ZFC}+V=L)$ implies $V=L$.  Also it is shown that if ${\rm LZFC}+Loc({\rm ZFC}+{\rm CH})+Loc({\rm ZFC}+\neg {\rm CH})$ is consistent, then   Powerset is false, while  the consistency of ${\rm LZFC}+ Loc({\rm ZFC}+{\rm CH})+Loc({\rm ZFC}+\neg {\rm CH})$ follows from the consistency of ZFC+``there is a natural model of ZFC''. Finally we show that for any definable set $c$ and definable ordinals $\alpha,\beta$, the theory
${\rm LZFC}+Loc({\rm ZFC}+|c|=\omega_\alpha)+ Loc({\rm ZFC}+|c|=\omega_\beta) +\ \mbox{\rm ``$c$ exists''}$ is inconsistent.

In section 7 we consider (in ZFC) a question  that arises as a
result of dealing  exclusively with transitive models. We can dub
it ``standard compactness'' problem, since it is like  ordinary
compactness except that the  models allowed  are  (standard)
transitive ones only.  Given a set $\Sigma$ of sentences of a
finitary language extending the language of set theory,  such that
$|\Sigma|=\kappa$ and every subset of $\Sigma$ of cardinality
$<\kappa$ has a transitive model, does $\Sigma$ have a transitive
model? If the answer is yes we  call $\kappa$ {\em standard
compact}. We show (in ZFC): (a) $\omega$ is not standard compact,
(b) every weakly compact cardinal is standard compact,  and (c) if
$\lambda>\omega$ is strongly compact, then every $\kappa\geq
\lambda$ such that $\kappa^{<\kappa}=\kappa$ is standard compact.

\section{A localized variant of ZFC.}
$V$ is the universe of sets.  The membership relation between
entities of $V$ is denoted by $\in$. Let ${\cal L}=\{\epsilon\}$
be the language of set theory. Since $\epsilon$ is going to be
interpreted only by $\in$  we shall identify $\epsilon$ with $\in$
and write for simplicity ${\cal L}=\{\in\}$.

$\Pi_n$, $\Sigma_n$ denote the usual classes of
formulas in  the  L\'{e}vy hierarchy (with
$\Pi_0=\Sigma_0$ being the class of bounded formulas).
If S is a set theory, $\Sigma_n^{\rm S}$ and  $\Pi_n^{\rm S}$ are the classes of formulas provably equivalent  in S to  a $\Sigma_n$ and  $\Pi_n$ formula, respectively. Also  $\Delta_n^{ \rm S}$ is the class of properties $\phi$ which are provably equivalent in S both to a $\Pi_n$ and a $\Sigma_n$ formula,  i.e.,  there is a $\Sigma_n$ formula $\phi_1$ and a $\Pi_n$  formula $\phi_2$  such
that ${\rm S}\vdash \phi\leftrightarrow \phi_1\leftrightarrow
\phi_2$.

Lower case letters $a,b, x, y, u, v$ denote sets. Upper case letters
$A,B,M,N,X,Y$ denote either sets or (proper) classes, depending on
the context. For example throughout the letters $M,N$ always
denote transitive {\em sets} which are models of ZFC.

If $\phi$ is a formula of ${\cal L}$ and $u$ is a set,  $\phi^{u}$
denotes the bounded formula resulting from $\phi$ if we replace
each unbounded quantifier $\forall x$, $\exists x$ of $\phi$ with
$\forall x\in u $, $\exists x\in u$, respectively. As usual
writing $\phi$ we mean that $(V,\in)\models \phi$. So $\phi^{u}$
is equivalent to $(u,\in)\models \phi$.

The following  localistic substitute of ZFC will be  the main
axiom of our system LZFC defined below:
$$(Loc({\rm ZFC)}) \hspace{.4\columnwidth minus .5\columnwidth}
 \forall x\exists y (x\in y \
\wedge \ Tr(y) \ \wedge \ (y,\in)\models {\rm ZFC}),
\hspace{.5\columnwidth minus .5\columnwidth} \llap{}$$ where
$Tr(y)$ denotes the formula ``$y$ is transitive'' and
$(y,\in)\models {\rm ZFC}$ abbreviates the formula   $\forall
\phi\in {\rm ZFC} ((y,\in)\models \phi)$. $Loc({\rm ZFC)}$ says
that the class of transitive models of  ZFC is an {\em unbounded}
(or {\em cofinal}) subclass of $V$ with respect to $\in$, and
hence with respect to $\subseteq$ (because of transitivity).

However, the relation ``$(y,\in)\models \phi$'', as well as the
set  ZFC, as a set of formulas,  cannot be defined without some
elementary notions and facts from a body of absolute set theoretic
truths that  we call  Basic Set Theory and  denote by BST. This
is similar to Elementary Set Theory, EST, of  \cite[p. 39]{En04},
except that BST contains in addition Cartesian Product, while the axioms
of Foundation and Choice are not included because
they can be deduced from $Loc({\rm ZFC})$ (see below).
So we take  BST to consist  of the following axioms:

(Emptyset) \ \ $\exists x(x=\emptyset)$,

(Ext) \ \ $\forall x\forall y[\forall z(z\in x\leftrightarrow z\in
y) \rightarrow x=y]$

(Pair) \ \ $\forall x\forall y\exists z (z=\{x,y\})$

(Union) \ \ $\forall x\exists y (y=\bigcup x)$

(Cartesian Product) \ \ $\forall x\forall y\exists z(z=x\times
y)$. [The predicates, ``pair'', ``function'' etc,  are $\Delta_0$ and can be
defined as in  \cite[p. 14]{Ba75}.]

(Infinity) \ \ $\exists x[\emptyset\in x \wedge \forall y\in
x\exists z\in x(z=y\cup\{y\})]$

$(\Delta_0$-Separation) \  \ $\forall \bar{z}\forall a\exists
b\forall y(y\in b\leftrightarrow y\in a \wedge \phi(y,\bar{z}))$,
\\ for every $\Delta_0$ formula $\phi$ not containing $b$ free.

\begin{Lem} \label{L:Peano}
In {\rm BST}: (i)  $\omega$ exists and the axioms of Peano
arithmetic {\rm (PA)} can be proven to hold in $\omega$ endowed
with the usual operations. Thus ${\rm PA}\subseteq {\rm BST}$.
(ii) The set of formulas $Fml({\cal L})$ is definable, $V_\omega$
exists and the relation ``$(x,\in)\models\phi(\bar{a})$''  is
definable.
\end{Lem}

{\em Proof.} (i) By Infinity, let $a$ be an inductive set. We can
define $\omega$ (using $\Delta_0$-Separation) as the set of
ordinals $x\in a$ such that for every  $y\leq x$ and $y\neq 0$,
$y$ is a successor ordinal. We can see that this set is the least
inductive set (details are left to the reader). The minimality of
$\omega$ as inductive set amounts  to the fact that $\omega$
satisfies complete induction. The operations $',+,\cdot$ on it are
defined as usual and  the axioms of PA are shown in BST to be true
with respect to $\omega$.

(ii) By Cartesian Product, for every set $a$ and $n\in\omega$,
$a^n=\{(x_0,\ldots,x_{n-1}):x_i\in a\}$ is a set.   Formulas of
${\cal L}$ are defined inductively as  triples of integers, e.g.
$\lceil v_i=v_j\rceil=(0,i,j)$,   $\lceil v_i\in
v_j\rceil=(1,i,j)$, etc, as in \cite[p. 90]{Dr74}.  The set
$Fml({\cal L})$ of formulas of ${\cal L}$ is a recursive (hence
$\Delta_1^{\rm ZFC}$ definable subset of $\omega$). So is also
${\rm ZFC}\subset Fml({\cal L})$.

The set $a^n$  can be identified also with the set of functions
$f$ such that $dom(f)=n$ and $rng(f)\subseteq a$. Using this
identification we can define $V_\omega$ as in \cite[p.81]{Dr74} by
a $\Delta_1^{\rm ZFC}$ definition. Finally the
relation``$(x,\in)\models \phi$'' is also  $\Delta_1^{\rm ZFC}$
definable by the help of $V_\omega$ (see \cite[p. 91]{Dr74} for
details). \telos

\vskip 0.2in

\begin{Rem} \label{R:cart}
{\em Without the axiom Cartesian Product of BST,  to prove that
cartesian products of sets are sets one  would need something like
$\Delta_0$-Collection (or $\Delta_0$-Replacement) (see \cite[prop.
3.2]{Ba75}). This is a rather strong axiom, while existence of
cartesian products is quite elementary. We do not know if LZFC proves $\Delta_0$-Collection (see Propositions \ref{P:equivalent} and  \ref{P:nice} below).}
\end{Rem}

Having fixed the  definitions of $Fml({\cal L})$, ZFC and
$(x,\in)\models \phi$, we can now consider the axiom $Loc({\rm
ZFC})$ given above and set
$${\rm LZFC}={\rm BST}+Loc({\rm ZFC}).$$
For simplicity  henceforth we shall write $x\models \phi$ instead
of $(x,\in)\models \phi$. Sometimes we drop also the predicate
$Tr(x)$ if implicitly understood, so $Loc({\rm ZFC})$ is usually
written $\forall x\exists y(x\in y \ \wedge \ y\models{\rm ZFC})$.

First let us note, as already mentioned above, that the  axioms of Choice and Foundation
are deduced from $Loc({\rm ZFC})$.

\begin{Lem} \label{L:rendered}
$Loc({\rm ZFC})$ implies the axioms of {\rm Choice} and {\rm Foundation}.
\end{Lem}

{\em Proof.} Let $x\neq \emptyset$ be a set such that for every $y\in x$, $y\neq\emptyset$. By  $Loc({\rm ZFC})$, there is a transitive model $M$ of ZFC such that $x\in M$. Then in $M$ $x$ has a choice function and also has a $\in$-least member. \telos

\vskip 0.2in

Given a tuple of sets $\bar x=(x_1,\ldots,x_n)$ let $\bar x\in y$
abbreviate the formula $x_1\in y\wedge\cdots \wedge x_n\in y$.

\begin{Lem} \label{L:singleax}

(i) ${\rm LZFC}\vdash \forall  \bar x \ \exists y(\bar x\in y
\wedge  Tr(y) \wedge y\models{\rm ZFC})$.

(ii)  Let $\Pi_2(\rm ZFC)$ be the set of $\Pi_2$ consequences of
{\rm ZFC}. Then  $\Pi_2({\rm ZFC})\subseteq {\rm LZFC}$.
\end{Lem}

{\em Proof.} (i) Let $\bar x=(x_1,\ldots,x_n)$. Given any
$a_1,\ldots,a_n$,  $\{a_1,\ldots,a_n\}$ exists in BST.  So by
$Loc({\rm ZFC})$ there is a transitive model $b$ such that
$\{a_1,\ldots,a_n\}\in b$. Then $\{a_1,\ldots,a_n\}\subset b$ and
$b\models{\rm ZFC}$.

(ii) Let $\phi\in \Pi_2(\rm ZFC)$. $\phi$ has the form $\forall
\bar{x}\ \exists\bar{y}\psi(\bar{x},\bar{y})$, where $\bar{x}$ is
an $n$-tuple of variables,  $\bar{y}$ is an $m$-tuple of variables
and  $\psi$ is bounded. Let us work in  ${\rm LZFC}$. Pick any $n$-tuple of sets $\bar{a}$. It suffices to show
that there is a $\bar{y}$ such that $\psi(\bar{a},\bar{y})$. By (i)
above  there is a transitive  model $b\models{\rm ZFC}$ such that
$\bar{a}\in b$. Since $\phi$ is a consequence of ZFC, $b\models
\phi$, or $b\models\exists\bar{y}\psi(\bar{a},\bar{y})$. Hence
$\exists\bar{y} \psi(\bar{a},\bar{y})$ since $\psi$ is $\Delta_0$.
Thus ${\rm LZFC}\vdash \phi$. \telos

\begin{Rem} \label{R1}
{\em  In contrast to ZFC, LZFC should not in general allow axioms
with unbounded quantifiers, since its  truths are ``local'', and
so  the variables must range in some set-model. However certain
$\Pi_2$ statements expressing  elementary, indisputable facts
(like e.g. $\forall x,y \exists z(z=\{x,y\})$), cannot but be
accepted, despite the occurrence of two alternating unbounded
quantifiers. This is the case with the axioms of BST. All of them
are $\Pi_2$ sentences, as one can easily check by inspecting the
formulations given above.}
\end{Rem}

\begin{Rem} \label{R2}
{\em  The axioms of BST are necessary  only to make possible the
strict}  formulation {\em of $Loc({\rm ZFC})$. Otherwise, that is,
if we assume that $Loc({\rm ZFC})$} is {\em sensible, by assuming
for example that the notions  ``formula'' and ``$x\models\phi$''
are primitive, then we can easily prove lemma \ref{L:singleax} by
working in   $Loc({\rm ZFC})$+Pair+Emptyset rather than LZFC.
Since all axioms of BST are $\Pi_2$ consequences of ZFC, it
follows from \ref{L:singleax} (ii) that from $Loc({\rm
ZFC})$+Pair+Emptyset we can recover the rest of the axioms of BST.
}
\end{Rem}

\begin{Rem} \label{R3}
{\em Throughout  we are going to make heavy use of the well-known
fact that every $\Delta_1^{\rm ZFC}$ (and hence every
$\Sigma_0^{\rm ZFC}$) formula is absolute for transitive models of
ZFC. However a word of caution is needed here.  $\Delta_1^{\rm
ZFC}$ formulas are absolute between transitive models of ZFC and
the universe, {\em when  we work in } ZFC (and this is done most
of the time), i.e., when  $V$ is supposed to satisfy ZFC. If
$V\not\models{\rm ZFC}$ absoluteness of $\Delta_1^{\rm ZFC}$
formulas  is no longer guaranteed. For instance let $\phi$ be
$\Sigma_1$ and $\psi$ be $\Pi_1$ and ${\rm ZFC}\vdash
\phi\leftrightarrow\psi$. Then $M\models \phi\leftrightarrow\psi$
for any model of ZFC. But if $V\not\models{\rm ZFC}$ we cannot
infer that $V\models \phi\leftrightarrow\psi$, so we cannot infer
absoluteness of $\phi$ and $\psi$.  In our case $V$ satisfies LZFC
rather than ZFC, so this observation is in order. However, if
${\rm S}$ is a set theory such that $V\models {\rm S}$ and ${\rm
S}\vdash \phi\leftrightarrow\psi$ whenever ${\rm ZFC}\vdash
\phi\leftrightarrow\psi$, for $\phi,\psi$ as above, then
$\phi,\psi$ are still absolute between $V$ and the models of ZFC.
The next lemma says that this is the case for ${\rm S=LZFC}$.}
\end{Rem}

\begin{Lem} \label{L:nobother}
If $\phi\in \Sigma_1$ and $\psi\in \Pi_1$ and ${\rm ZFC}\vdash
\phi\leftrightarrow \psi$, then ${\rm LZFC}\vdash
\phi\leftrightarrow \psi$.  Consequently for any transitive
$M\models {\rm ZFC}$, any $\Delta^{\rm ZFC}_1$ formula
$\phi(\bar{x})$ and any $\bar{a}\in M$,
$\phi(\bar{a})\leftrightarrow M\models \phi(\bar{a})$.
\end{Lem}

{\em Proof.} Let $\phi\in \Sigma_1$ and $\psi\in \Pi_1$ and ${\rm
ZFC}\vdash \phi\leftrightarrow \psi$. Then $\phi\leftrightarrow
\psi$ belongs to $\Pi_2({\rm ZFC})$, so the claim follows from
lemma  \ref{L:singleax} (ii).  \telos

\vskip 0.2in

For brevity  we express the fact established in lemma
\ref{L:nobother} by saying that every  $\Delta_1^{\rm ZFC}$
formula of ${\cal L}$ is also $\Delta_1^{\rm LZFC}$.

It follows from lemma \ref{L:nobother} that  every set defined by
a $\Delta^{\rm ZFC}_1$ formula inside any transitive model $M$ of
ZFC with parameters in $M$ is the same as when defined in  ${\rm LZFC}$. We often express this by saying that this
set {\em exists in} LZFC, in the sense that its  definition in  LZFC
does not provide a proper class.  In particular this is the case
with sets defined inductively by some positive inductive operator
$\Gamma_\phi$ in any transitive model, for some $\Sigma_1$ formula
$\phi$.

\begin{Rem} \label{R:LZ}
{\em Let us remark at this point, for later use, that the sentence
$$Loc({\rm ZFC})=\forall x\exists y(x\in y\wedge
y\models{\rm ZFC})$$  is itself $\Pi_2^{\rm ZFC}$, since ``$\phi$ is a
formula'', ``$\phi\in {\rm ZFC}$'' and ``$x\models\phi$'' are
$\Delta_1^{\rm ZFC}$. Moreover, by lemma \ref{L:nobother},
$Loc({\rm ZFC})$ is also  $\Pi_2^{\rm LZFC}$. }
\end{Rem}

Ordinals are defined  in LZFC as usual (transitive sets linearly ordered, and hence well-ordered, by $\in$). Lower case Greek letters $\alpha,\beta,\ldots$  denote ordinals. We
often write $\alpha<\beta$ instead of $\alpha\in \beta$. We denote the
class of all ordinals  by $On$.  $(On,\in)$ is well-ordered, but we must be careful
with the meaning of this assertion.  $(On,\in)$ is well-ordered means that every {\em subset} of $On$ has a least element, as a consequence of Foundation. Things however may be different for {\em subclasses} of $On$. If $X=\{\alpha:\phi(\alpha)\}$ is a subclass of $On$, then there is no way to ensure that $X$ has a least element. The usual argument that amounts  to pick an $\alpha\in X$  and then take the trace $\alpha\cap X$ of $X$ on $\alpha$ does not work in LZFC since, in absence of full separation,   $\alpha\cap X$ need not be  a set. It works only for $\Delta_0$-classes (i.e., classes defined by $\Delta_0$-formulas).  So let us denote  by ${\rm Found}_{On}$ the scheme ``every subclass of $On$ has a least element''. Namely:

\vskip 0.1in

$({\rm Found}_{On})  \ \  \exists \alpha\in On \ \phi(\alpha) \rightarrow \exists
\alpha\in On[\phi(\alpha) \wedge \forall \beta<\alpha \neg\phi(\beta)]$.

\vskip 0.1in

${\rm Found}_{On}$ is apparently a weak form of  the full $\in$-induction scheme ${\rm Found}^*$ which says that ``every class has an $\in$-least element'':

\vskip 0.1in

$({\rm Found}^*) \ \  \exists x\phi(x) \rightarrow \exists
x[\phi(x) \wedge \forall y\in x\neg\phi(y)]$

\vskip 0.1in
However we shall see   below (Lemma \ref{L:classorder}) that ${\rm Found}_{On}$ and ${\rm Found}^*$ are  in fact equivalent over LZFC.

A  remarkable  situation where  ${\rm Found}_{On}$ is involved is the following. Let us call  sets $x,y$ {\em equinumerous} and write $x\sim y$ if there is a bijection $f:x\rightarrow y$. Also let us write $x\precsim y$ if there is an injection $f:x\rightarrow y$, and $x\precnsim y$ if there is an injection $f:x\rightarrow y$, but $x\not\sim y$.  Given any $x$, let $Ord(x)=\{\alpha\in On:x\sim \alpha\}$. By Choice, for every $x$, $Ord(x)\neq \emptyset$. However the formula $x\sim \alpha$ is $\Sigma_1$, hence, since $\Sigma_1$-Separation is not available in LZFC,  we cannot ensure that $Ord(x)$ has a least element. The least element of $Ord(x)$, if it existed, would be  the (absolute) cardinality of $x$, what we usually denote $|x|$. It follows that in LZFC alone, without ${\rm Found}_{On}$ (or ${\rm Found}^*$ according to the previous lemma), absolute cardinalities of sets cannot be defined.\footnote{Consequently the notation $|x|$ will not be used when $x$ is a set of LZFC.  Sometimes this notation is employed without actual reference to existent cardinalities as sets. For example, the notation $|x|=|y|$ is another way  to say   $x\sim y$, while  $|x|<|y|$ means just $x\precnsim y$.}  This is rather in accordance with the spirit of LZFC, whose primary motivation was to challenge  the existence of absolute infinite cardinalities and powersets. So further discussion on this issue is provided in section 3.

Given a model $M$, let $Def(M)$ denote the collection  of its
first-order definable subsets, i.e., $$Def(M)=\{X\subseteq
M:(\exists \phi(x,\bar{y})\in \omega)(\exists \bar{b}\in
M)[M\models \forall x(x\in X\leftrightarrow \phi(x,\bar{b}))]\}.$$
The definition  is absolute so $Def(M)$ exists in $V$. Further, if
${\cal X}$ is a subset of ${\cal P}(M)$, then $Def(M,{\cal X})$
denotes the collection of subsets of $M$ second-order definable in
$(M,{\cal X})$.

The Ramified Analytical hierarchy over $M$ is the collection
$RA(M)=\bigcup_{\alpha\in On}RA_\alpha(M)$, where

$RA_0(M)=Def(M)$,

$RA_{\alpha+1}(M)=Def(M,RA_\alpha(M))$,

$RA_\alpha(M)=\bigcup_{\beta<\alpha}RA_\beta(M)$.

\begin{Lem} \label{L:list}
The following facts are provable in ${\rm LZFC}$ and are absolute
with respect to transitive models of {\rm ZFC}:

(i) $\exists x(x=V_\omega)$ (the set of hereditarily finite sets
exists).

(ii) $\forall x\exists y (TC(x)=y)$ (every set has a transitive
closure).

(iii) $\forall x\exists \alpha\in On (rank(x)=\alpha)$, where
$rank(x)=\sup\{rank(y)+1:y\in x\}$.

(iv) $\forall \alpha \in On \ \exists x(x=L_\alpha)$, where
$L_\alpha$ is the $\alpha$-th level of the ordinary constructible
hierarchy.

(v) $\Delta_1$-Separation.

(vi) For every model $M$, $(M,Def(M))$ as well as $(M,RA(M))$
exist and are models of the theories of classes {\rm GBC}
(G\"{o}del-Bernays) and {\rm KM} (Kelley-Morse), respectively.

\end{Lem}

{\em Proof.} All objects involved in the clauses (i)-(vi) above
have  $\Delta^{\rm ZFC}_1$ definitions,  therefore   $\Delta^{\rm
LZFC}_1$ definitions by \ref{L:nobother}, and hence they have absolute definitions inside any transitive model of ZFC containing the appropriate parameters. For instance to show existence
of $L_\alpha$, take a transitive model $M\models{\rm ZFC}$ such
that $\alpha\in M$, and construct in $M$ the levels $L_{\beta}$,
$\beta\leq \alpha$, of $L$. \telos

\vskip 0.2in

\begin{Lem} \label{L:classorder}
 ${\rm Found}^*$ and  ${\rm Found}_{On}$ are equivalent over ${\rm LZFC}$.
\end{Lem}

{\em Proof.} Since the ordering $<$ on $On$ coincides with $\in$, obviously ${\rm Found}^*$ implies  ${\rm Found}_{On}$. Conversely, suppose
${\rm Found}_{On}$ holds and let $\exists x\phi(x)$ be true. Consider the subclass of $On$  $$R_\phi=\{\alpha\in On:\exists x (\phi(x) \ \wedge \ rank(x)=\alpha)\}.$$
By Lemma \ref{L:list} (iii),  every set in LZFC has a rank, hence  $R_\phi\neq \emptyset$.  By ${\rm Found}_{On}$, $R_\phi$ has a least element $\alpha_0$. Thus $\exists x (\phi(x) \wedge rank(x)=\alpha_0)$ is true. Pick such a  $x$. Then  $\forall y\in x \neg \phi(y)$. \telos

\vskip 0.2in

In view of the  non-derivability of ${\rm Found}_{On}$ in LZFC  we have the following important consequence.

\begin{Rem} \label{R:noinduction}
{\em The familiar transfinite induction along $On$ is not available in {\rm LZFC}, except for $\Delta_1$ subclasses of $On$. }
\end{Rem}

Because of   \ref{L:list} (iii),  we can  define (non-inductively!) for every
$\alpha\in On$ the class
$$V_\alpha=\{x:rank(x)<\alpha\}.$$ $V_\alpha$, $\alpha\in On$, are
the {\em layers} of the universe, since $V_\alpha\subset V_\beta$
for $\alpha<\beta$ and $V=\bigcup_\alpha V_\alpha$. Except
$V_\alpha$ for $\alpha\leq \omega$, $V_\alpha$ in general need not
be sets. However it is straightforward that  the relativization of
$V_\alpha$'s to any transitive model $M$ of ZFC generates the
usual cumulative hierarchy  of $M$.

\begin{Lem} \label{L:cumulative}
Let $M$ be a transitive model of ${\rm ZFC}$. Then for every
$\alpha\in On^M$, $V_\alpha^M=M_\alpha=M\cap V_\alpha$.
\end{Lem}

Let also
$$L=\bigcup_{\alpha\in On}L_\alpha$$
be the class of constructible sets. In contrast to $V_\alpha$,
each $L_\alpha$ is a set.\footnote{However one should not expect
LZFC to prove  that $L$ is an inner model of ZFC, since for that
one would need the Replacement Axiom. After all such a requirement
would not comply with the localistic spirit of LZFC, according to
which only {\em set} models of ZFC make sense.} So the picture of
the universe of ${\rm LZFC}$ is roughly that of Figure 1.

\vskip 0.5in

\begin{center}
\setlength{\unitlength}{1mm}
\begin{picture}(30,30)(-20,-30)
\linethickness{1pt} \thinlines

\put(-35,1){\line(0,-1){25}} \put(35,1){\line(0,-1){25}}
\put(-35,-24){\line(1,0){70}} \put(-10,-30){\line(2,-3){10}}
\put(10,-30){\line(-2,-3){10}}

\put(-10,-30){\line(-1,4){8}}

\put(10,-30){\line(1,4){8}}

\put(0,1){\line(0,-1){45}} \put(-35,-12){\line(1,0){70}}
\put(-35,-17){\line(1,0){70}} \put(-35,-6){\line(1,0){70}}
\put(13,-31){\line(0,-1){4}} \put(13,-40){\line(0,-1){5}}
\put(36,-12){\line(1,0){5}} \put(41,-13){\line(0,-1){15}}
\put(41,-33){\line(0,-1){12}} \put(-10,-30){\line(1,0){20}}
\put(11,-30){\line(1,0){2}} \put(35,-26){\line(0,-1){7}}
\put(35,-40){\line(0,-1){5}} \put(35,-25){\line(-5,-1){24}}
 \put(-35,-25){\line(5,-1){24}}

\put(-2,-10){\makebox(0,0)[c]{$\alpha$}}
\put(-6,-23){\makebox(2,2)[c]{\small $\omega+1$}}
\put(-2,-28){\makebox(0,0)[c]{\small $\omega$}}
\put(-4,0){\makebox(0,0)[c]{$On$}}
\put(14,-38){\makebox(1,1)[c]{$V_\omega$}}
\put(0,-48){\makebox(1,1)[c]{$\emptyset$}}
\put(41,-32){\makebox(2,2)[c]{$V_\alpha$}}
\put(-32,0){\makebox(0,0)[c]{$V$}}
\put(-14,0){\makebox(0,0)[c]{$L$}}
\put(34,-35){\makebox(0,0)[c]{$V_{\omega+1}$}}

\end{picture}
\end{center}

\vskip 0.8in

\begin{center}
Figure 1
\end{center}

\begin{Rem} \label{R:pow}
{\em The picture of Figure 1 suggests that the levels $V_\alpha$
for $\alpha>\omega$ are all proper classes. This however need not
be always true and some $V_\alpha$ may be sets in some cases. Of
course if $V_\alpha$ is a proper class, so is every $V_\beta$ for
$\beta>\alpha$.  LZFC simply
does not give any information  about the status of the Powerset
axiom and, by so doing, is generally compatible with ZFC (see
\ref{P:PC} below), although its intended interpretation points to
the opposite direction. In order to   refute the Powerset
axiom, we need localization principles stronger than $Loc({\rm
ZFC})$, of the form  $Loc({\rm ZFC}+\phi)$ or $\{Loc({\rm
ZFC}+\phi):\phi\in\Gamma\}$.    See section 6.}
\end{Rem}

Beside the Powerset axiom,  the axiom scheme of Collection/Replacement is also questionable when referred to $V$. In general for a set of formulas $\Gamma$ we have the scheme:

\vskip0.1in

$(\Gamma$-Collection) \  $\forall \bar{z}\forall a\exists
b[\forall x\in a \exists y\phi(x,y,\bar{z}) \rightarrow \forall
x\in a\exists y\in b \ \phi(x,y,\bar{z})]$,  for every formula
$\phi\in\Gamma$ not containing $b$ free.

\vskip0.1in

$\Gamma$-Replacement is  weaker than $\Gamma$-Collection, so we
consider only the latter.

\vskip 0.1in

Also the following  scheme of $\Pi_2$-Reflection is of interest here:

\vskip 0.1in

$(\Pi_2$-Reflection) \  $\phi\rightarrow \ \forall x\exists y[x\in
y \wedge Tr(y) \wedge  \phi^{y}]$,  for every $\Pi_2$ sentence
$\phi$.

\vskip 0.1in

(Clearly, working in LZFC  we may use in the above scheme
$\Pi_2^{\rm LZFC}$ sentences instead of just $\Pi_2$.)

\begin{Prop} \label{P:equivalent}
(i) $\Delta_0{\mbox{-} \rm Collection}$ and  $\Sigma_1{\mbox{-} \rm Collection}$ are equivalent over ${\rm LZFC}$.

(ii) ${\rm LZFC}+\Pi_2{\mbox{-}\rm Reflection}\vdash
\Sigma_1{\mbox{-}\rm Collection}$.
\end{Prop}

{\em Proof.} (i) One direction is trivial. It suffices to show that  $\Delta_0{\mbox{-} \rm Collection}$ implies  $\Sigma_1{\mbox{-} \rm Collection}$ over LZFC. Let $\psi(x,y):=\exists \bar{z}\phi(x,y,\bar{z})$ be a $\Sigma_1$-formula, and  let $\forall x\in a\exists y\psi(x,y)$ be true in LZFC. Then $\forall x\in a\exists y\exists \bar{z}\phi(x,y,\bar{z})$. Let $n$ be the length of the tuple $\bar{z}$. Using  pairing and the $\Delta_0$ functions $(u)_0,\ldots,(u)_n$, for an $(n+1)$-tuple $u$,  such that $u=((u)_0,\ldots,(u)_n)$,
$\forall x\in a\exists y\exists \bar{z}\phi(x,y,\bar{z})$ is written $\forall x\in a\exists u\phi(x,(u)_0,(u)_1,\ldots,(u)_n)$. Since $\phi(x,(u)_0,(u)_1,\ldots,(u)_n)$ is (an abbreviation of) a $\Delta_0$ formula, by $\Delta_0$-Collection there is a $b$ such that $\forall x\in a\exists u\in b \ \phi(x,(u)_0,(u)_1,\ldots,(u)_n)$. If $c=TC(b)$, then
$\forall x\in a\exists y\in c \exists \bar{z}\phi(x,y,\bar{z})$, i.e., $\forall x\in a\exists y\in c \ \psi(x,y)$.

(ii) We work in ${\rm LZFC}+\Pi_2$-Reflection. Let
$\phi(x,y,\bar{z})$ be a $\Sigma_1$ formula, $a$, $\bar{c}$ be
sets and let $\forall x\in a \exists y\phi(x,y,\bar{c})$ hold
true. We have to show that there is $b$ such that
$$\forall x\in a\exists y\in b \
\phi(x,y,\bar{c}).$$  Since $\phi$ is $\Sigma_1$,  $\forall x\in a
\exists y\phi(x,y,\bar{c})$ is a $\Pi_2$ formula. By
$\Pi_2$-Reflection there is a transitive $b$ such that
$a\cup\{c_1,\ldots,c_n\}\in b$ and $(\forall x\in a \exists
y\phi(x,y,\bar{c}))^{b}$, or $\forall x\in a \exists y\in b \
\phi(x,y,\bar{c})^{b}$. Since $\phi$ is $\Sigma_1$,
$\phi^{b}$ implies $\phi$, so  $\forall x\in a \exists y\in b \
\phi(x,y,\bar{c})$. \telos

\vskip 0.2in

Let ${\rm TM}({\rm LZFC})$ denote the principle ``there is a transitive
model of LZFC''.

\begin{Prop} \label{P:nice}
${\rm LZFC+\Pi_2}$-{\rm Reflection} \ $\vdash {\rm TM}({\rm LZFC})$.
Consequently ${\rm LZFC+\Pi_2}$-{\rm Reflection} \ $\vdash
Con({\rm LZFC})$. Therefore if {\rm LZFC} is consistent, then
${\rm LZFC}\not\vdash \Pi_2$-{\rm Reflection}.
\end{Prop}

{\em Proof.} We work in ${\rm LZFC+\Pi_2}$-{\rm Reflection}. By
Remark  \ref{R:LZ},  the axiom $Loc({\rm ZFC})$ of LZFC is a true
$\Pi_2$ sentence, hence $\Pi_2$-Reflection applies to
$Loc({\rm ZFC})$.  Consider the conjunction $\Phi=Loc({\rm
ZFC})\wedge {\rm Pair}$. Clearly $\Phi$ is $\Pi_2$, and by
assumption it holds in $V$, so by $\Pi_2$-{\rm Reflection} there
is a (nonempty) transitive set $b$ such that $\Phi^b=(Loc({\rm
ZFC}))^b\wedge {\rm Pair}^b$ is the case.  It suffices to show
that $b\models {\rm LZFC}$. Already $b\models Loc({\rm ZFC})
\wedge {\rm Pair}$, so it remains to show that $b$ satisfies the
rest of the axioms of ${\rm BST}$. Emptyset and Extensionality are
obvious in view of the transitivity of $b$. For Union, let $x\in
b$. Then $x\in M\in b$ for some model $M$, so $\bigcup x\in M\in
b$. For Cartesian Product, given any  $x,y\in b$, $\{x,y\}\in b$
by Pair, so there is, by $Loc({\rm ZFC})$, a model $M\in b$ such
that $\{x,y\}\in M$. Then $x,y\in M$, hence $x\times y\in M\in b$.
Similarly for Infinity. It remains to
verify $\Delta_0$-Separation. Let $c\in b$ and $\phi(x,\bar{a})$
be $\Delta_0$, with $\bar{a}\in b$.   Let $X=\{x\in
c:b\models\phi(x,\bar{a})\}$. We have to show that $X\in b$. By
Pair and $Loc({\rm ZFC})$ there is a model $M\in b$ of ZFC such
that $c,\bar{a}\in M$. Then clearly $X=\{x\in M:M\models x\in c
\wedge \phi(x,\bar{a})\}$. Therefore $X\in M$ and hence $X\in b$.
The proof that ${\rm LZFC+\Pi_2}$-{\rm Reflection} \ $\vdash
TM({\rm LZFC})$ is complete. So ${\rm LZFC+\Pi_2}$-{\rm
Reflection} \ $\vdash Con({\rm LZFC})$. Since by lemma
\ref{L:Peano} ${\rm PA}\subseteq {\rm LZFC}$, G\"{o}del's
incompleteness implies that ${\rm LZFC}\not\vdash \Pi_2$-{\rm
Reflection}.  \telos

\vskip 0.2in

It is open whether   ${\rm LZFC}$ proves $\Delta_0{\mbox{-}\rm Collection}$. Also it is open  whether the converse of \ref{P:equivalent} (ii) above is true, i.e., whether ${\rm LZFC}+\Sigma_1{\mbox{-}\rm Collection}$ proves $\Pi_2{\mbox{-}\rm Reflection}$. (If it does, then, in view of  \ref{P:equivalent} (i) and \ref{P:nice}, ${\rm LZFC}\not\vdash\Delta_0{\mbox{-}\rm Collection}$).

As a byproduct of the proof of the last proposition we have the
following simple fact that gives a sufficient condition in order
for a set to be a model of LZFC. A transitive set $(a,\in)$ is
said to be {\em directed} if it is upward directed as a poset,
i.e., if  for all $x,y\in a$ there is a $z\in a$ such that $x,y\in
z$.

\begin{Cor} \label{C:light}
Let $a$ be a transitive set which is the union of the transitive
models of ${\rm ZFC}$ contained in it, that is, $a=\bigcup\{x\in
a:x\models{\rm ZFC}\}$. If $a$ satisfies also  Pair, then
$a\models{\rm LZFC}$. In particular, if $(a,\in)$ is a  directed
set of models of {\rm ZFC}, such that $\cup a=a$, then
$a\models{\rm LZFC}$.
\end{Cor}

In the preceding result we can even replace models of ZFC with
models of LZFC. Namely the following holds.

\begin{Lem} \label{L:directed}
Let $(a,\in)$ be a  directed set of models of {\rm LZFC}, such
that $\cup a=a$. Then $a\models{\rm LZFC}$.
\end{Lem}

{\em Proof.} By directedness $a$ satisfies {\em Pair.} So   it
suffices to show that $a\models Loc({\rm ZFC})$. Let $x\in a$.
Then there is $b\in a$ such that $x\in b$ and $b\models Loc({\rm
ZFC})$. Therefore $b\models \exists y (x\in y \wedge y\models{\rm
ZFC})$. But then $a\models (\exists y (x\in y \wedge y\models{\rm
ZFC}))^b$, hence $a\models \exists y (x\in y \wedge y\models{\rm
ZFC})$, or $a\models  Loc({\rm ZFC})$. \telos

\vskip 0.2in

Clearly if ZFC and LZFC are consistent theories, then ${\rm
ZFC}\not\subseteq {\rm LZFC}$ and ${\rm LZFC}\not\subseteq {\rm
ZFC}$.  Of the other set theories of the literature, close to the
BST part of  ${\rm LZFC}$ is Kripke-Platek set theory with
infinity (KP + Infinity) (see \cite{Ba75}, where rather the system
KPU=KP+ urelements is considered). This is the system of axioms:
$${\rm KP}=\{{\rm Empty, Ext,  Pair, Union, {Found}^*},
\Delta_0{\rm\mbox{-}Separation}, \Delta_0{\rm\mbox{-}Collection}\},$$ where
${\rm Found}^*$ is the  scheme of $\in$-induction we already saw  above to be  equivalent to ${\rm Found}_{On}$ (see Lemma \ref{L:classorder} and before) and does not seem to follow from LZFC. It follows that
$${\rm KP+ Infinity}-\{\Delta_0{\rm\mbox{-}Collection}, {\rm Found}^*\}\subset {\rm LZFC}.$$
In connection with Remark \ref{R:noinduction}, let us cite here the reasonable  extensions of LZFC in which induction is valid.

\begin{Lem} \label{L:indvalid}
$${\rm LZFC}+{\rm Found}_{On}\subseteq {\rm LZFC}+{\rm Separation}\subseteq {\rm LZFC}+{\rm Replacement}\subseteq {\rm LZFC}+{\rm Collection}.$$
\end{Lem}

{\em Proof.} The first inclusion follows from the discussion after Remark \ref{R:LZ}. The third inclusion is obvious.  Concerning the inclusion ${\rm LZFC}+{\rm Separation}\subseteq {\rm LZFC}+{\rm Replacement}$, the proof is no different from the familiar one that is used in ZFC. \telos

\vskip 0.2in

The systems  LZFC + Separation and LZFC + Replacement, apart from the fact that they restore transfinite induction, seem to be interesting in themselves extensions of LZFC.

A few further existence results for LZFC are given below.

\begin{Lem} \label{L:Skolem}
(i) The L\"{o}wenheim-Skolem theorem is provable in ${\rm LZFC}$.
Namely, for every  first-order language ${\cal L}$, every ${\cal
L}$-structure ${\cal A}=(A,\ldots)$ and every  $S\subseteq A$ such
that $S\precsim {\cal L}$, there is an  ${\cal L}$-structure
${\cal B}=(B,\ldots)$ such that $S\subseteq B$, $B\precsim {\cal L}$
and ${\cal B}\preceq {\cal A}$.

(ii) The Mostowski's isomorphism theorem  is provable in ${\rm
LZFC}$. Namely if $x$ is a set and $E$ is a binary relation on $x$
such that (a) $E$ is well-founded and (b) $(x,E)\models {\rm
Ext}$, then there is a (unique) transitive set $y$ such that
$(x,E)\cong (y,\in)$.

(iii) The Completeness Theorem is provable in ${\rm LZFC}$.
\end{Lem}

{\em Proof.} All three theorems, when formalized, are $\Pi_2$
sentences provable in ZFC, so the claim follows from lemma
\ref{L:singleax}. \telos

\begin{Lem} \label{L:CA}
 ${\rm ACA}\subset {\rm LZFC}$.
\end{Lem}

{\em Proof.}  First-order Peano axioms, when transcribed into
${\cal L}=\{\in\}$, become $\Delta_0$ sentences, since all
quantifiers are restricted to $\omega$. The induction axiom
$$\forall X[(0\in X \wedge \forall n(n\in X \rightarrow n+1\in X))
\rightarrow  \forall n(n\in X)]$$ becomes a $\Pi_1$ sentence,
since $\forall X$ becomes $\forall x\subseteq \omega$. The
arithmetic comprehension axiom is
$$\forall \bar{X} \ \exists Y\forall n(n\in Y\leftrightarrow
\phi(n,\bar{X})),$$ where $\phi$ has no set quantifiers. In ${\cal
L}=\{\in\}$ it  becomes $$\forall \bar{x} \ \exists y\forall z\in
\omega(z\in y\leftrightarrow \psi(z,\bar{x})),$$ where $\psi$ now
is bounded, hence  a $\Pi_2$ sentence provable in ZFC. Thus, in
view of \ref{L:singleax} (ii), both the induction axiom and the
comprehension scheme of ACA are provable in LZFC, hence  ${\rm
ACA}\subset {\rm LZFC}$. \telos

\begin{Lem} \label{L:byforcing}
(i) Let $M$ be a countable transitive  model of {\rm ZFC} and let
$B\in M$ be a Boolean algebra. Then  it is provable in ${\rm
LZFC}$ that there are $M$-generic filters $G\subseteq B$.

(ii) For every $M$ and generic $G$ as above the generic extension
$M[G]$ exists in ${\rm LZFC}$.
\end{Lem}

{\em Proof.} (i) Given a countable $M$ and the algebra  $B\in M$,
an $M$-generic filter $G\subseteq B$ is constructed by Choice  as
usual.

(ii) Given $M$, $B$ an $G$ as above, $M[G]$ is constructed by two
inductive definitions: One that provides  the set $M^B$ of
$B$-names over $M$, and another that leads from $M^B$ and $G$ to
the $G$-interpretations of $M^B$, $I_G''M^B=M[G]$. Both
definitions are inductive and absolute. So carrying them out
inside any model $N$ such that $M,B,G\in N$, is the same as
carrying them out in $V$.  \telos

\vskip 0.3in

{\bf Consistency.} What about the truth and  consistency of ${\rm
LZFC}$? Let ${\rm IC}$ be the  axiom ``there exists a strongly
inaccessible cardinal'',  ${\rm IC}^\infty$ be the  axiom
``there is a proper class of strongly inaccessible cardinals'' and  NM be the
axiom ``there is a natural model of ZFC'' (i.e., of the form $V_\alpha$). It is well known that the implications ${\rm IC}^{\infty}\rightarrow {\rm IC}\rightarrow {\rm NM}$  are strict over ZFC.

\begin{Prop} \label{P:PC}
(i) ${\rm LZFC}\subset {\rm ZFC}+{\rm IC}^\infty$.

(ii) ${\rm ZFC}+{\rm IC}\vdash Con({\rm ZFC}+{\rm LZFC})$.

(iii) ${\rm ZFC+}{\rm NM}\vdash Con({\rm LZFC+\mbox{\rm ``Every set is countable''}})$.
\end{Prop}

{\em Proof.} (i) Work in ${\rm ZFC}+{\rm IC}^\infty$. It suffices to
prove that $Loc({\rm ZFC})$ holds.  Then every set $a$ belongs to
some $V_\kappa$, where $\kappa$ is strongly inaccessible. Since
every such $V_\kappa$ is a transitive model of ZFC, it follows
that $\forall x\exists y(x\in y \wedge \ y\models {\rm ZFC})$.

(ii) Let $\kappa$ be an inaccessible in the ZFC universe. Then $V_\kappa\models {\rm ZFC}+Loc({\rm ZFC})$. Indeed, obviously $V_\kappa\models {\rm ZFC}$. It is well known (see \cite{MV59}, or \cite[Ex. 12.12]{Je03}) that $\{\alpha\in V_\kappa: (V_\alpha,\in)\prec (V_\kappa,\in)\}$ is closed unbounded in $\kappa$. Hence $\forall x\in V_\kappa \exists y(x\in y \wedge y\models {\rm ZFC})$. Thus $V_\kappa\models Loc({\rm ZFC})$.

(iii)  Let $V_\kappa$ be a natural model of ZFC. It is well-known that $\kappa$ is sufficiently large so that $H(\omega_1)\in V_\kappa$. $H(\omega_1)$ is the required model. Indeed, let  $x\in H(\omega_1)$. Then   $x\in V_\kappa$ and by L\"{o}wenheim-Skolem there is a countable model  $N\prec V_\kappa$ such that $x\in N$. If  $N'$ is the  Mostowski collapse of $N$, then   $N'$ is a  transitive model that contains $x$ and belongs to $H(\omega_1)$. Therefore,  $H(\omega_1)\models Loc({\rm ZFC})$. Moreover $H(\omega_1)\models$ ``Every set is countable''. \telos

\vskip 0.2in

It follows from \ref{P:PC} (ii) that the consistency strength of
${\rm LZFC}$ is no greater than that of ${\rm ZFC}+{\rm NM}$. Also by  \ref{P:PC} (iii), the consistency of ${\rm ZFC}+Loc({\rm ZFC})$ is no greater than that of ${\rm ZFC}+{\rm IC}$. Moreover,  ${\rm ZFC}+Loc({\rm ZFC})$
is a good mild substitute of ${\rm ZFC}+{\rm IC}^\infty$. It's worth
mentioning that ${\rm IC}^\infty$ is equivalent to what in category
theory is called ``the axiom of universes'', the origin of which
goes back to Grothendieck. Roughly a ``Grothendieck universe'' is
a transitive set closed under pairing, powerset and replacement.
The axiom of universes says that every set belongs to a
Grothendieck universe. It is likely that most or all of what the
category theorists prove by the help of the axiom of universes,
can be proved within ${\rm ZFC}+Loc({\rm ZFC})$.

\section{Cardinals and powersets in ${\rm LZFC}$}
Typically,  we may keep talking about cardinals in LZFC, much the same way as we do in ZFC, but without expecting to prove the familiar ZFC results,  due to the lack of  Powerset, Replacement and also transfinite induction (Remark \ref{R:noinduction}). The landscape of LZFC is hazy as far as absolute infinite cardinalities are concerned,  and pitfalls are lurking everywhere for the visitor accustomed to ZFC.

We can define cardinals as usual. An ordinal $\alpha$ is a said to be a {\em cardinal} (in the sense of $V$) if it is an initial ordinal, i.e., if there is no $\beta\in On$ such that $\beta<\alpha$ and $\beta\sim \alpha$. For instance $\omega$ is a cardinal. In fact $\omega$ may be the only infinite cardinal (as it follows from  Proposition  \ref{P:PC} (ii)).
$\omega_1$ is the class of countable ordinals,  i.e., $$\omega_1=\{\alpha\in On:\alpha\precsim \omega\}.$$ In general this  is a (proper) class. A class $X=\{x:\phi(x)\}$ is said to {\em exist,} if it is a set. So if $\omega_1$ exists,  it is a cardinal. If $\omega_1$ is a proper class,  can we infer that $\omega_1=On?$ Actually not. Because  $\omega_1$ is an initial segment of $On$, but in order to draw a contradiction from  $On-\omega_1\neq \emptyset$, the latter class should have a least element, which we cannot guarantee. If $\beta\in On-\omega_1$, then $\omega_1$ would be a  proper subclass of $\beta$. If $M$ is a model of ZFC containing $\beta$, then $\omega_1\subseteq \beta\subseteq M$, but $\omega_1^M\subsetneqq \omega_1$, i.e., $\omega_1^M\in \omega_1$, otherwise $\omega_1=\omega_1^M$ and   $\omega_1$ would have to be a set.

If $\omega_1$ exists, then we set
$\omega_2=\{\alpha\in On:\alpha\precsim \omega_1\}$,
and  similar remarks apply to this class. If $\omega_1$ exists, then by $Loc({\rm ZFC})$ there is a model $M$ of ZFC such that $\omega_1\in M$. $\omega_1$ is clearly a cardinal in $M$ but not necessarily an absolute one  with respect to $M$. It may be the case that $\omega_1^M<\omega_1$, and hence $\omega_1=\omega_\alpha^M$, for some $\alpha>1$. But even if $\omega_1^M=\omega_1$, $\omega_2^M$ (which is a set) need not be absolute, and $\omega_2$ may be a proper class. In general, the class $X=\{\alpha\in On: \omega_\alpha \ \mbox{exists}\}$ is defined, but we know neither  whether  $X=On$ nor whether $On-X$ has a least element.

Analogous comments hold  about the power-class ${\cal P}(\omega)$ and the class
$$H(\omega_1)=\{x:TC(x)\precsim \omega\}$$ of hereditarily countable sets.
If  ${\cal P}(\omega)$ exists, then  ${\cal P}(\omega)$ belongs to a model $M$ and, obviously,  ${\cal P}(\omega)^M={\cal P}(\omega)$. Since however $M$ need not be a natural model, it is possible that   $({\cal P}^2(\omega))^M\neq {\cal P}^2(\omega)$ and, moreover, ${\cal P}^2(\omega)$ be a proper class. Again for the class $Y=\{\alpha:{\cal P}^\alpha(\omega) \ \mbox{exists}\}$ we can say neither whether $Y=On$, nor whether $On-Y$ has a least element. Also,  if ${\cal P}^n(\omega)$ exists for every $n\in \omega$, we cannot conclude that ${\cal P}^\omega(\omega)$ exists, since  Replacement is missing.

Concerning $H(\omega_1)$, it is well-known that in ZFC we can code its elements  by elements of ${\cal P}(\omega)$, constructing thus an embedding $f:H(\omega_1)\rightarrow {\cal P}(\omega)$. This is done by induction on the rank of the elements of $H(\omega_1)$ which  goes up to $\omega_1$. So this embedding cannot be carried out in  LZFC.

The above  uncertainties  about absolute infinite cardinalities  seem to fit to the spirit of LZFC. They prompt one to deal exclusively with models and let aside  absolute uncountable infinities. However the  uncertainties are settled as soon as we augment  LZFC with  Separation, which restores transfinite induction (see Lemma \ref{L:indvalid}).

At this point   I would like  to address an ambiguity (that occurs also in the ZFC environment), concerning the meaning of the symbols  $\omega_\alpha$.   $\omega_\alpha$ is allowed to  denote alternatively (depending on the context) either  an {\em object,} i.e., a specific  ordinal, or a {\em property,}  the property of being the $\alpha$-th infinite cardinal. The ambiguity arises from  the interplay of the two meanings within models of  {\rm ZFC}. For instance if $M\models {\rm ZFC}$,  $\beta\in On\cap M$ and we write $M\models \beta=\omega_\alpha$, we refer to $\omega_\alpha$  as a property, namely, the property ``$\beta$ is  the $\alpha$-th infinite cardinal number'' (in the sense of $M$). The last assertion  is alternatively  denoted $\beta=\omega_\alpha^M$. Similarly, in the expression $M\models |x|=\omega_\alpha$, $\omega_\alpha$ is construed as a property. Now assume that $\omega_\alpha$ is a set. By $Loc({\rm ZFC})$  there is a model $M$ such that  $\omega_\alpha\in M$. If for some $x\in M$ we write $M\models x\sim \omega_\alpha$, then we refer to $\omega_\alpha$ as an object which is involved in a property that is true in $M$.  On the other hand, $\omega_\alpha$ is  still a cardinal in $M$, but it need not preserve also its size, i.e., we may have $\omega_\alpha=\omega_\beta^M$ for some $\beta>\alpha$. According to the usage of $\omega_\beta$ as a property, the latter is written equivalently  $M\models \omega_\alpha=\omega_\beta$, which seems to be absurd. The absurdity is simply  due to the ambiguity of the symbols $\omega_\alpha,\omega_\beta$: In the formula $M\models \omega_\alpha=\omega_\beta$, $\omega_\alpha$ is construed as an object, while $\omega_\beta$ is construed as a property. The situation is no different in ZFC. Simply the (set) models we deal with  there are, mostly, either {\em countable}, hence they do not contain real  uncountable cardinals, or  {\em natural,} in which all powersets and cardinals are absolute.  The problematic situation is exactly when $\omega_\alpha$ is uncountable, $\omega_\alpha\in M$ and  $\omega_\alpha^M\neq \omega_\alpha$.

We can raise the  ambiguity if we avoid using the symbols $\omega_\alpha$ as properties and employ instead a predicate $Card(\alpha,x)$ for the property ``$x$ is the $\alpha$-th infinite cardinal number''.  The predicate  $Card(\alpha,x)$ is defined as follows. Let $$Card(x):= x\in On \wedge  \forall \beta<x (\beta\not\sim x)$$ be the property ``$x$ is a cardinal''. Then the formula $Card(\alpha,x)$ is defined by the following clauses:
$$\left\{\begin{array}{l}
           Card(0,x):=[x=\omega] \\
            Card(\alpha+1,x):= [Card(x) \wedge \forall y (Card(\alpha,y)\rightarrow y\precnsim x \ \wedge \\ \hspace{1.3in} \forall z(Card(z) \rightarrow z\precsim y \vee x\precsim z))] \\
             Card(\alpha,x):=[(\forall \beta<\alpha \forall y(Card(\beta,y)\rightarrow y\precnsim x) \ \wedge   \\
            \hspace{1.1in} \forall z(Card(z) \rightarrow x\precsim z \vee \exists \gamma<\alpha\exists u(Card(\gamma,u) \wedge z\precsim u))],\\
            \hspace{1.1in} \mbox{for $\alpha$ limit}.
            \end{array} \right.
$$
Note that $Card(\alpha,x)$ is intended to be used inside models of ZFC, so the induction on $\alpha$ needed to verify $M\models Card(\alpha,x)$ is legitimate. Using the predicate $Card(\alpha,x)$, we  write $M\models Card(\alpha,\beta)$ instead of $M\models \beta=\omega_\alpha$. If $\omega_\alpha\in M$ and $\omega_\alpha$ happens to be the $\beta$-th cardinal of $M$, we express it by writing $M\models Card(\beta,\omega_\alpha)$  instead of the puzzling $M\models \omega_\alpha=\omega_\beta$. This way the ambiguity is removed.

Below we shall keep using the notation $\beta=\omega_\alpha^M$ as an abbreviation of  $M\models Card(\alpha,\beta)$. Also $M\models |x|=\omega_\alpha$ will be  an abbreviation of
$$M\models \exists \beta (x\sim \beta \wedge Card(\alpha,\beta)).$$ If $\omega_\alpha$ exists and $M$ is a model such that $\omega_\alpha\in M$, we say that $\omega_\alpha$ is {\em absolute in $M$} if $\omega_\alpha^M=\omega_\alpha$, i.e., if $M\models Card(\alpha,\omega_\alpha)$.
The following is easy to verify.

\begin{Lem} \label{L:Galois}
{\rm (LZFC)}  If $M,N$ are models of {\rm ZFC} such that $M\subseteq N$, $x\in M$, $\alpha\in M$, and  $M\models |x|=\omega_\alpha$, then  $N\models |x|\leq \omega_\alpha$.
\end{Lem}

{\em Proof.} We just argue as usual inside the model $N$.  \telos

\vskip 0.2in

In general, if $n\in \omega$ and  $\omega_n$ exists, we set   $H(\omega_{n+1})=\{x:TC(x)\precsim\omega_n\}$ and $\omega_{n+1}=\{\alpha\in On:\alpha\precsim \omega_n\}$.

\begin{Lem} \label{L:settle}
In {\rm LZFC}, for all $n\in \omega$, the following hold.

(i) If $H(\omega_{n+1})$ exists, then so do  ${\cal P}(\omega_n)$ and $\omega_{n+1}$. In particular, if $M$ is a model of {\rm ZFC} such that $H(\omega_{n+1})\in M$, then  ${\cal P}(\omega_n)^M={\cal P}(\omega_n)$ and $\omega_{n+1}^M=\omega_{n+1}$.

(ii) If  ${\cal P}(\omega_n)$ exists and ${\cal P}(\omega_n)\in M$ then $\omega_{n+1}^M=\omega_{n+1}$ and $H(\omega_{n+1})^M=H(\omega_{n+1})$.

(iii) Suppose  ${\cal P}^{n+1}(\omega)$ exists and ${\cal P}^{n+1}(\omega)\in M$.   Then ${\cal P}(\omega_n)\in M$, hence $\omega_{n+1}^M=\omega_{n+1}$. Also  ${\cal P}^{n+1}(V_\omega)^M=V_{\omega+n+1}^M=V_{\omega+n+1}$.

Moreover in  {\rm LZFC+Separation}, the above claims are proved for every $\alpha\in On$.  Namely:

(iv) If  $H(\omega_{\alpha+1})\in M$, then  ${\cal P}(\omega_\alpha)^M={\cal P}(\omega_\alpha)$ and $\omega_{\alpha+1}^M=\omega_{\alpha+1}$.

(v) If ${\cal P}^\alpha(\omega)\in M$, then $\omega_\alpha^M=\omega_\alpha$ and $V_{\omega+\alpha}^M=V_{\omega+\alpha}={\cal P}^\alpha(V_\omega)$.
\end{Lem}

{\em Proof.} For clarity and simplicity we show  clauses (i) and (ii) for $n=0$ and clause (iii) for $n=1$. The inductive steps are straightforward and left to the reader.

(i) Suppose  $H(\omega_1)$ is a set and $M$  is a model such that $H(\omega_1)\in M$. Then   ${\cal P}(\omega)\subseteq H(\omega_1)\subseteq M$, therefore ${\cal P}(\omega)= {\cal P}(\omega)^M$.  Also, $\omega_1^M=\{\alpha\in On\cap M:M\models \alpha\precsim \omega\}$. Hence $\omega_1^M\subseteq \omega_1$. For the converse,  let $\alpha\in \omega_1$ be an infinite ordinal. Then there is a bijection  $f:\alpha\rightarrow \omega$. Clearly $f\in H(\omega_1)$, and hence $f\in M$. Since $\alpha=dom(f)$, $\alpha\in M$, therefore $\alpha\in \omega_1^M$. So $\omega_1^M=\omega_1$.

(ii) Suppose ${\cal P}(\omega)$ exists and let ${\cal P}(\omega)\in M$. We show first that $\omega_1^M=\omega_1$. As we saw above, $\omega_1^M \subseteq \omega_1$. To show the converse, pick some infinite $\alpha\in \omega_1$. It suffices to show that $\alpha\in M$ and $M\models \alpha\sim \omega$.  Now  there is (in $V$) a bijection $f:\omega\rightarrow \alpha$. Let  $$R=\{\langle m,n\rangle\in \omega\times\omega:f(m)\in f(n)\}.$$
$\omega\times\omega$ is a set and the defining property of $R$ is $\Delta_0$, so by $\Delta_0$-Separation, $R$ is a set too. Moreover $R$  is a well-ordering of $\omega$ and $R\in {\cal P}(\omega\times\omega)$. Since ${\cal P}(\omega)\in M$, also ${\cal P}(\omega\times\omega)\in M$. Hence $R\in M$ and $M\models ``(\omega,R) \ \mbox{is a well-ordering''}$. So the order type of $(\omega,R)$ exists in $M$. But  this order-type is $\alpha$, i.e.,   $\alpha\in M$ and $M\models (\alpha,\in)\cong (\omega,R)$. Therefore $M\models \alpha\sim \omega$.

We come to the second claim of this clause, and let ${\cal P}(\omega)\in M$. We have to show that  $H(\omega_1)^M=H(\omega_1)$, or $H(\omega_1)\subseteq H(\omega_1)^M$.\footnote{The proof of this implication was provided by the referee.} Let $x\in H(\omega_1)$, and let $f:TC(x)\rightarrow \omega$ be a bijection. Let $N$ be a model of ZFC such that $\{{\cal P}(\omega),f\}\subset  N$. In $N$ we can define as usual a coding $g:H(\omega_1)^N\rightarrow {\cal P}(\omega)^N={\cal P}(\omega)$. Now  the pair $\langle x,f\rangle$ is an element of $H(\omega_1)^N$ and it is coded by $g(\langle x,f\rangle)\in {\cal P}(\omega)$. But since ${\cal P}(\omega)\in M$, $g(\langle x,f\rangle)$ is in $M$ and  from $g(\langle x,f\rangle)$ we can fully restore $\langle x,f\rangle$, i.e., $\langle x,f\rangle\in M$. Thus $x\in  H(\omega_1)^M$.

(iii) We show the claim  for $n=1$. Let ${\cal P}^2(\omega)\in M$.  Then ${\cal P}(\omega)\in M$, and hence $\omega_1^M=\omega_1$, by (ii). Every $\alpha\in \omega_1$ is coded by some well-ordering $R\in{\cal P}(\omega\times \omega)$ of $\omega$, as we saw in (ii). Hence every $x\subseteq \omega_1$ is coded by some element of ${\cal P}^2(\omega \times\omega)$, or equivalently, of ${\cal P}^2(\omega)$. So ${\cal P}(\omega_1)$ is (coded by) a subset of ${\cal P}^2(\omega)$. This means that ${\cal P}(\omega_1)\in M$  and,  by (ii), $\omega_2^M=\omega_2$.  The other claim also follows easily.

(iv) and (v) need induction on $\alpha$.
Here we cannot work in any particular model of LZFC,  so the induction must be carried out  in $V$. This explains the use of Separation.
\telos

\vskip 0.2in

Concerning the converse of  the claims (i)-(iii)  above, some of them  can be shown to be false (assuming the consistency of some basic theory). For instance  it is consistent relative to ZFC +LZFC that in LZFC $\omega_1$ exists, while ${\cal P}(\omega)$ is a proper class. Indeed, if  ZFC +LZFC is consistent, then so is ZFC +LZFC +${\cal P}(\omega)\sim \omega_2$. If $K$ is a model of the last theory, then  $H(\omega_2)^K$ is a model of LZFC + ``$\omega_1$ exists'' + ``${\cal P}(\omega)$ does not exist''.

\section{Mahlo models}
Transitive models of ${\rm ZFC}$ bear obvious analogies with
inaccessible cardinals. Roughly a transitive $M\models {\rm ZFC}$
is a ``first-order counterpart'' of an inaccessible cardinal,
since both  are transitive sets closed under the same basic closure conditions. These closure conditions are related with the two most powerful axioms of ZFC, Replacement and Powerset. First, a (strongly)  inaccessible cardinal $\kappa$ is closed under all functions $f$, in the sense that for every $\alpha\in \kappa$, $f''\alpha$ is bounded in $\kappa$. The corresponding property of a model $M$ is that, in view of Replacement, for every $x\in M$, $f''x\in M$, provided $f$ is {\em first-order definable} in $M$. (That is what we mean by saying that $M$ is  a first-order counterpart of an inaccessible cardinal). Second, for every cardinal $\lambda<\kappa$, $2^\lambda<\kappa$, and this obviously corresponds to the truth of Powerset in $M$, i.e., the fact that for every $x\in M$, ${\cal P}^M(x)\in M$.\footnote{Even in  ZFC, the existence of a transitive model of
ZFC can be thought as a weak large cardinal axiom, in view of the
non-reversible implications
$${\rm IC}\rightarrow {\rm NM} \rightarrow {\rm TM} \rightarrow Cns({\rm ZFC}),$$
where

IC ``There is an inaccessible cardinal'',

NM: ``There is a natural (i.e., of the form
$V_\alpha$) model of ZFC'',

TM: ``There is a transitive model of ZFC'',

$Cns({\rm ZFC})$: ``ZFC is consistent''. }

Consequently,   a transitive $M$ such that $M\models {\rm ZFC}+Loc({\rm ZFC})$ is the analogue of a ``quasi 1-Mahlo'' cardinal in the following
sense: $M\models Loc({\rm ZFC})$ says that every $x\in M$ belongs
to a $y\in M$ such that  $y\models{\rm ZFC}$. That is, the set of
transitive models contained  in $M$ form an unbounded (= cofinal)
subclass of $M$ under $\in$ (and $\subseteq$). This is just the
property of being   $1$-Mahlo cardinal, except that ``unbounded''
should be replaced by ``stationary''. So $M$ is  ``quasi 2-Mahlo''
if $M\models {\rm ZFC}+Loc({\rm ZFC}+ Loc({\rm ZFC}))$, and so
on.\footnote{Note that  the operator $Loc$ can be applied not only
to ZFC, but to any set theory S in the obvious way. Namely
$Loc({\rm S}):=\forall x\exists y(x\in y \wedge y\models {\rm S})$.
In order however for the latter to make sense,  S must be a  {\em
definable} set of axioms in a language ${\cal L}'\supseteq {\cal
L}$. If S is defined by $\phi$, then $Loc({\rm S})$ is the ${\cal
L}'$-sentence $\forall x\exists y(x\in y\wedge \forall
z(\phi(z)\rightarrow y\models z))$.} Stationarity, however, is  a
relative notion: It depends on what closed unbounded sets (clubs)
are available. Absoluteness is obtained only if one is confined to
the collection of  {\em definable} clubs  and stationary subsets
of a model $M$. Before coming to the definition of stationary
subsets of models, let us define inductively the iterated
localization principles  $Loc_n({\rm ZFC})$, for $n\in\omega$, as
follows:

$ Loc_0({\rm ZFC})=Loc({\rm ZFC})$,

$Loc_{n+1}({\rm ZFC})= Loc({\rm ZFC}+Loc_n({\rm ZFC)})$. \\
It is easy to check that  for every $n\in\omega$, the sentence
$Loc_n({\rm ZFC)}$ is $\Pi_2$. Moreover inductively we can see
that
\begin{equation} \label{E:stronger}
 Loc_{n+1}({\rm ZFC}) \ \rightarrow  \ Loc_n({\rm ZFC}).
\end{equation}

\begin{Rem} \label{R:wherestop}
{\em Can we continue the definition of $Loc_\alpha({\rm ZFC})$ for
$\alpha\geq \omega$?  The definition can be carried out at
least along the constructive ordinals in a way analogous to that
used in \cite{Fe62} for the definition of transfinite progressions
of theories using the consistency operator:   ${\rm T}_0={\rm T}$, ${\rm
T}_{\alpha+1}={\rm T}_\alpha+Con({\rm T}_\alpha)$, ${\rm
T}_\alpha=\bigcup_{\beta<\alpha}{\rm T}_\beta$. In ZFC one may
also define $Loc_\alpha({\rm ZFC})$ by using ordinals
$\beta<\alpha$ as parameters. For example suppose that
$Loc_\beta({\rm ZFC})$, for $\beta<\alpha$, have  been defined, so
that the mapping $\beta\mapsto Loc_\beta({\rm ZFC})$ is definable.
Then, by Replacement, $\{ Loc_\beta({\rm ZFC}):\beta<\alpha\}$ is
a definable set, so, in view of footnote 7,  we can set
$Loc_\alpha({\rm ZFC})=Loc({\rm ZFC}\cup\{Loc_\beta({\rm
ZFC}):\beta<\alpha\})$. But in {\rm LZFC} Replacement is not
available, so $\{Loc_\beta({\rm ZFC}):\beta<\alpha\}$ need not be a
set, and therefore iteration of $Loc$ cannot go beyond
constructive ordinals.}
\end{Rem}

Recall that  a cardinal $\kappa$ is said to be Mahlo if the set of
inaccessibles below $\kappa$  is stationary in $\kappa$.  Since
the axioms $Loc_n({\rm ZFC)}$ involve only the unboundedness of
the class of models, just Mahloness of $\kappa$ suffices in order
for $V_\kappa$ to satisfy $Loc_n({\rm ZFC)}$.

\begin{Prop} \label{P:Mahlo}
{\rm (ZFC)} Let  $\kappa$ be a Mahlo cardinal. Then
$V_\kappa\models Loc_n({\rm ZFC})$ for all $n\in \omega$.
\end{Prop}

{\em Proof.} Let us define inductively for $n\in\omega$, that a
cardinal  $\kappa$ is $n$-{\em unbounded} if it is inaccessible
and for every $m<n$, the $m$-unbounded cardinals are unbounded in
$\kappa$.

{\em Claim 1.} If $\kappa$ is Mahlo, then $\kappa$ is
$n$-unbounded for all $n\in \omega$.

{\em Proof.} By induction on $n$. Trivially  $\kappa$ is $0$- and
$1$-unbounded. Suppose $\kappa$  is $n$-unbounded for $n\geq 1$.
Then  the $(n-1)$- unbounded cardinals are cofinal in $\kappa$.
Let $\alpha<\kappa$. Let   $S$ be the set of limit points of
$(n-1)$-unbounded above $\alpha$.  It is easy to check that
$S$ is a club. So, since $\kappa$ is Mahlo, $S$ contains an
inaccessible $\beta$. This $\beta$ is also a limit of
$(n-1)$-unbounded cardinals, so it is an $n$-unbounded and lies  above
$\alpha$. This means that the $n$-unbounded cardinals are cofinal
in $\kappa$. Hence $\kappa$ is $(n+1)$-unbounded.

{\em Claim 2.}  If $\kappa$ is $(n+1)$-unbounded, then
$V_\kappa\models Loc_n({\rm ZFC})$.

{\em Proof.} By induction on $n$.  Let $\kappa$ be $1$-unbounded.
Then the set $S\subset\kappa$ of inaccessibles below $\kappa$ is
unbounded in $\kappa$. For every $\lambda\in S$,
$V_\lambda\models{\rm ZFC}$. Therefore $V_\kappa$ satisfies
$\forall x\exists y(x\in y \wedge y\models{\rm ZFC})$, i.e.,
$V_\kappa\models Loc_0({\rm ZFC})$.

We assume that the claim holds for $n+1$ and we show it for $n+2$.
Let $\kappa$ be $(n+2)$-unbounded. The set $S\subset \kappa$ of
$(n+1)$-unbounded cardinals   is unbounded in $\kappa$. By the
induction hypothesis, for every $\lambda\in S$, $V_\lambda\models
Loc_n({\rm ZFC})$. Therefore $V_\kappa$ satisfies $\forall
x\exists y(x\in y \wedge y\models{\rm ZFC}+Loc_n({\rm ZFC}))$. The
last sentence is $Loc({\rm ZFC}+Loc_n({\rm ZFC}))=Loc_{n+1}({\rm
ZFC})$.

Claims 1 and 2 yield the proof of the proposition. \telos

\vskip 0.2in

The iterated localization principles $Loc_n({\rm ZFC})$ are ``weak
Mahlo'' principles  intended to motivate the full Mahlo notion for
models considered below. The latter presumes the notion of club
and stationary set adapted here for that purpose. Unless otherwise
stated, the definitions below are given in LZFC.

\begin{Def} \label{D:club}
{\em Let  $M$ be a transitive model of ${\rm ZFC}$. A set $X\in
Def(M)$  is said to be} unbounded in $M$, {\em if $(\forall x\in
M)(\exists y\in X)( x\subseteq y)$. A $X\in Def(M)$   is said to
be} closed, {\em if
$$(\forall y\in M)(y\subseteq X \wedge (y,\subseteq) \ \mbox{is
a chain} \ \rightarrow \ \cup y\in X).$$ A $X\in Def(M)$  is said
to be a} club {\em  of $M$ if it is unbounded and closed. A $X\in
Def(M)$ is said to be} stationary in $M$ {\em if $X\cap Y\neq
\emptyset$ for every club $Y\in Def(M)$.}
\end{Def}

For a model   $M\models{\rm ZFC}$, a typical club of $M$ is the
set
$$\{M_\alpha:\alpha\in On\cap M\},$$ where $M_\alpha=V_\alpha^M$.
For every $M\models{\rm ZFC}$, let
$$Club(M)=\{x\in Def(M): x \ \mbox{is closed unbounded in } \
M\},$$
$$Stat(M)=\{x\in Def(M): x \ \mbox{is stationary in } \
M\}. $$ Since  $Def(M)$ is absolute, it follows that $Club(M)$ and
$Stat(M)$  are absolute too. It is easy to see that for every $M$,
$Club(M)$ is a proper subset of $Stat(M)$. For instance, if $X$ is
a club, $(y,\subseteq)$ is a chain of $X$ and we set $Y=X-\{\cup
y\}$, then $Y\in Stat(M) \backslash Club(M)$.

Given a transitive $M\models{\rm ZFC}$ and any  unbounded $X\in
Def(M)$, let $F^M_X:On^M\rightarrow On^M$ be defined as follows:
$$F^M_X(\alpha)=\mbox{\rm least}\{\beta:(\exists x\in X)
(M_\alpha\subseteq x\subseteq M_\beta)\}.$$ Clearly $F^M_X\in
Def(M)$. $F^M_X$  is said to be  the {\em associated function} to
$X$ with respect to $M$. We write simply $F_X$ instead of $F_X^M$
whenever $M$ is understood. It follows from the definition that
\begin{equation} \label{E:central}
(\forall \alpha\in On^M)(\exists x\in X)(M_\alpha\subseteq
x\subseteq M_{F_X(\alpha)}).
\end{equation}

\begin{Lem} \label{L:increase}
For every $M$ and every definable unbounded $X\subseteq M$, (a)
$F_X$ is nondecreasing, i.e., for all $\alpha<\beta\in M$,
$F_X(\alpha)\leq F_X(\beta)$. (b) For every $\alpha\in M$,
$\alpha\leq F_X(\alpha)$.
\end{Lem}

{\em Proof.} (a) Let $\alpha<\beta$. Then $\exists x\in X
(M_\beta\subseteq x\subseteq M_{F_X(\beta)})$. Since
$M_\alpha\subseteq M_\beta$, we have $\exists x\in X
(M_\alpha\subseteq x\subseteq M_{F_X(\beta)})$. Since
$F_X(\alpha)$ is the least $\gamma$ such that $\exists x\in X
(M_\alpha\subseteq x\subseteq M_\gamma)$, it follows that
$F_X(\alpha)\leq F_X(\beta)$. (b) Just  note that, by definition,
$M_\alpha\subseteq M_{F_X(\alpha)}$, therefore $\alpha\leq
F_X(\alpha)$. \telos

\vskip 0.2in

With the help of the function $F_X$ one can prove the following
closure properties of clubs. Since they are not going to be used
in the proof of the main result of the section, Proposition
\ref{P:sconnection}, we omit the proofs.

\begin{Lem} \label{L:meet}
(i) For any  $X_1,X_2\in Club(M)$,  $X_1\cap\ X_2\in Club(M)$.

(ii)  Let $X\in Def(M)$ be a set of pairs coding a family of clubs
of $M$. i.e., for every $x\in dom(X)$, $X_{(x)}=\{y:(x,y)\in X\}$
is a club. Then for every set $A\subseteq dom(X)$, $A\in M$,
$\bigcap_{x\in A}X_{(x)}\in Club(M)$.

(iii)  If $X_{(x)},x\in M$, is   an $M$-family of clubs of $M$,
then $\triangle_{x\in M}X_{(x)}$ is a club (where $\triangle_{x\in
M}X_{(x)}$ is the usual diagonal intersection of $X_{(x)}$). A
{\sl fortiori} $\triangle_{x\in S}X_{(x)}$ is a club for every
$S\in Def(M)$.
\end{Lem}

We come to the definition of $\alpha$-Mahlo models of ZFC.

\begin{Def} \label{D:Mahlo1}
{\em (LZFC)} $\alpha$-Mahlo models {\em of {\rm ZFC} are defined
inductively as follows:

(i) $x$ is} $0$-Mahlo {\em if $x$ is transitive and $x\models{\rm
ZFC}$.

(ii) $x$ is} $(\alpha+1)$-Mahlo, {\em  if $x$ is transitive,
$x\models{\rm ZFC}$ and $\{y\in x: (y,\in~) \ \mbox{is an
$\alpha$-Mahlo model}\}$  is a stationary subset of $x$.

(iii) For  $\alpha$ limit, $x$ is} $\alpha$-Mahlo {\em if it is
$\beta$-Mahlo for all $\beta<\alpha$.}
\end{Def}

The above definition of $\alpha$-Mahloness is  formalized by the formula $mahlo(\alpha,x)$ defined by the following clauses (we omit only transitivity of $x$  as implicitly understood):

\begin{equation} \label{E:informal}
\left\{\begin{array}{l}
           mahlo(0,x):=[x\models{\rm ZFC}] \\
            mahlo(\alpha+1,x):=[x\models{\rm ZFC} \wedge (\forall y\in
             Club(x))(\exists u\in y)(mahlo(\alpha,u))] \\
            mahlo(\alpha,x):=\forall \beta<\alpha \ mahlo(\beta,x), \ \mbox{for} \ \alpha \ \mbox{limit}.
            \end{array} \right.
\end{equation}
The lack of induction on $\alpha$ does not prevent $mahlo(\alpha,x)$ from having a truth value for all $\alpha$ and $x$. This is because $mahlo(\alpha,x)$ is absolute,  since $Club(x)$ is a $\Delta_1$ property. Hence the induction on $\alpha$ needed to verify $mahlo(\alpha,x)$ can be carried out inside any model $M$ containing $x$ and $\alpha$.

\begin{Lem} \label{L:second}
{\rm (LZFC)} For each $\alpha$, the sentence  $mahlo(\alpha,x)$ is
first-order and absolute for transitive models. That is, for every
transitive model  $M\models{\rm ZFC}$ such that $\alpha, x\in M$,
$mahlo(\alpha,x)$ iff $M\models mahlo(\alpha,x)$.
\end{Lem}

{\em Proof.} By an easy induction on $\alpha$, taking into account
that the right-hand sides of the clauses of (\ref{E:informal}) are
absolute. \telos

\vskip 0.2in

Note that  Mahloness alone (i.e., 1-Mahloness) implies the
iterated localization axiom $Loc_n({\rm ZFC})$.

\begin{Prop} \label{P:sameas}
{\rm (LZFC)} For every $n\in\omega$, if $M$ is Mahlo then
$M\models Loc_n({\rm ZFC})$.

\end{Prop}

{\em Proof.} The proof is similar to that of proposition
\ref{P:Mahlo} so it is omitted. \telos

\vskip 0.2in

Recall  that the clubs of a cardinal $\kappa$ are exactly the
ranges of normal (i.e., strictly increasing and continuous)
functions $f:\kappa\rightarrow \kappa$ (see e.g. \cite[p.
92]{Je03}). For every unbounded  $X\subseteq M$ (in particular for
every club), we defined above (see (\ref{E:central})) the
associated function $F_X:On^M\rightarrow On^M$, which is
nondecreasing rather than strictly increasing, and satisfies
$F_X(\alpha)\geq \alpha$. Such functions can  also be called
normal when they are continuous.\footnote{ If $f$ is simply
nondecreasing, i.e., $\alpha<\beta\rightarrow f(\alpha)\leq
f(\beta)$,  $rng(f)$ may  be bounded, which trivializes $f$. But
if  $rng(f)$ {\em is} unbounded, e.g. if $f(\alpha)\geq \alpha$,
then strictness of  monotonicity can be relaxed.  This is the case
with functions $F_X$. }  Obviously every  such function has
fixed points above any ordinal, as usual.  Using clubs $X$ such
that $F_X$ is normal, we can relate clubs of $M$ with clubs of
$\kappa$.

\begin{Def} \label{D:normalclub}
{\em Call a club  $X\subseteq M$} normal, {\em if the associated
function $F_X$ is normal.}
\end{Def}

Given $M\models {\rm ZFC}$, let $$U_M=\{M_\alpha:\alpha\in M\}$$
be the typical  club of $M$. For every $X\in Club(M)$, let us set
$$X^*=X\cap U_M.$$
By lemma \ref{L:meet} $X^*\in Club(M)$.

\begin{Lem} \label{L:normal}
For every $X\in Club(M)$, $X^*$ is a normal club.
\end{Lem}

{\em Proof.} Since, by \ref{L:increase} (b),  $F_{X^*}$ is already
nondecreasing, it suffices to show that $F_{X^*}$ is continuous,
i.e.,  for every limit $\alpha$,
$F_{X^*}(\alpha)=\sup\{F_{X^*}(\beta):\beta<\alpha\}$. Now the
elements of $X^*$ are sets $M_\beta$. Let $X^-=\{\beta\in
On^M:M_\beta\in X^*\}$. Then, by definition, for every $\beta$,
$$F_{X^*}(\beta)=\mbox{least}\{\gamma: (\exists x\in X^*)
M_\beta\subseteq x\subseteq M_\gamma\}=\mbox{least}\{\gamma\in
X^-:M_\beta\subseteq M_\gamma\}.$$ Therefore, for every
$\beta<\alpha$, $M_\beta\subseteq M_{F_{X^*}(\beta)}$ and
$M_{F_{X^*}(\beta)}\in X^*$. So
\begin{equation} \label{E:gamma}
\bigcup_{\beta<\alpha}M_\beta=M_\alpha\subseteq
\bigcup_{\beta<\alpha}M_{F_{X^*}(\beta)}=M_\gamma,
\end{equation}
where $\sup\{F_{X^*}(\beta):\beta<\alpha\}=\gamma$. But
$\{M_{F_{X^*}(\beta)}:\beta<\alpha\}\subseteq X^*$, and the chain
$\{M_{F_{X^*}(\beta)}:\beta<\alpha\}$ is in $M$. So, since $X^*$
is a club, $\bigcup_{\beta<\alpha}M_{F_{X^*}(\beta)}=M_\gamma\in
X^*$. Then  (\ref{E:gamma}) implies $F_{X^*}(\alpha)\leq \gamma$.
On the other hand, by monotonicity of $F_{X^*}$ (see
\ref{L:increase} (a)), $F_{X^*}(\alpha)\geq
\sup\{F_{X^*}(\beta):\beta<\alpha\}=\gamma$. So
$F_{X^*}(\alpha)=\gamma$ as required. \telos

\vskip 0.2in

In view of lemma \ref{L:normal}, a definable $Y\subseteq M$ is
stationary iff it meets all normal clubs of $M$ of the form $X^*$
for $X\in Club(M)$. For every model $M\models{\rm ZFC}$, let
$ht(M)$ (the height of $M$) be the supremum of the ordinals in
$M$, that is,  $ht(M)=M\cap On$.

Recall that

(i) $\kappa$ is  $0$-Mahlo  if it is strongly inaccessible.

(ii) $\kappa$ is  $(\alpha+1)$-Mahlo,  if  the  set of
$\alpha$-Mahlo cardinals below $\kappa$ is a stationary subset of
$\kappa$.

(iii) For limit $\alpha$,  $\kappa$ is  $\alpha$-Mahlo  if it is
$\beta$-Mahlo for all $\beta<\alpha$.

\begin{Prop} \label{P:sconnection}
(i) Let $M\models {\rm ZFC}$ with $ht(M)=\alpha$. If $X\in
Club(M)$, then $\{\beta<\alpha: M_\beta\in X^*\}$ is a club of
$\alpha$.

(ii) {\rm (ZFC)} If $\kappa$ is $\beta$-Mahlo, for $\beta<\kappa$,
then $V_\kappa$ is $\beta$-Mahlo.
\end{Prop}

{\em Proof.} (i)  Let $X\in Club(M)$ and  let
$X^-=\{\beta<\alpha:M_\beta\in X^*\}$. We have to show that $X^-$
is a club of $\alpha$. Let $\beta<\alpha$. It is clear that
$rng(F_{X^*})\subseteq X^-$. Since $F_{X^*}$ is normal, it has a
fixed point $\gamma>\beta$. Now $F_{X^*}(\gamma)=\gamma$ means
that $\gamma\in X^-$, so $X^-$ is unbounded. Further, let
$\{\beta_\xi:\xi<\delta\}$ be  an increasing sequence of $X^-$.
Then $\{M_{\beta_\xi}:\xi<\delta\}$ is an increasing sequence of
$X^*$. If $\beta=\sup\{\beta_\xi:\xi<\delta\}$, then
$M_\beta=\bigcup_{\xi<\delta}M_{\beta_\xi}$, and  $M_\beta\in
X^*$, by the closedness of $X^*$. Therefore  $\beta\in X^-$ and
$X^-$ is closed.

(ii) By induction on $\beta$. If $\kappa$ is $0$-Mahlo, then
$\kappa$ is strongly inaccessible, hence $V_\kappa\models{\rm
ZFC}$, and thus $V_\kappa$ is a $0$-Mahlo model according to
(\ref{E:informal}).

Suppose the claim holds for $\beta$ and let  $\kappa$ be
$(\beta+1)$-Mahlo. Then the set $Y=\{\lambda<\kappa:\lambda \
\mbox{is $\beta$-Mahlo}\}$ is stationary in $\kappa$. Let
$Y^+=\{V_\lambda:\lambda\in Y\}$. Both $Y$ and $Y^+$ are definable
in $V_\kappa$. By the induction hypothesis,  for every $x\in Y^+$,
$x$ is a $\beta$-Mahlo model. So  it suffices to show that $Y^+$
is stationary in $V_\kappa$, or, in view of \ref{L:normal}, that
it meets all  clubs $X^*$ for $X\in Club(V_\kappa)$. Let $X\in
Club(V_\kappa)$. Since $V_\kappa\cap On=\kappa$, by (i), the set
$X^-=\{\alpha<\kappa:V_\alpha\in X^*\}$ is a club of $\kappa$.
Therefore $Y\cap X^-\neq\emptyset$, hence  $Y^+\cap
X^*\neq\emptyset$.

If $\beta$ is limit then the claim follows immediately from the
definitions. \telos

\section{$\Pi_1^1$-Indescribable models}
The next question is whether models resembling higher large
cardinals can be reasonably defined. After Mahlo the next
candidate notion is that of a weakly compact model. However as is
well-known weakly compact cardinals have several equivalent
characterizations, through a partition property, a tree property,
a compactness property, $\Pi_1^1$-indescribability, etc (see for
example \cite{Je03}, \S 17). Although the  most intuitively
appealing characterization is  the partition property, the one
that  seems to fit better to our context is
$\Pi_1^1$-indescribability. Recall that a cardinal $\kappa$ is
$\Pi^n_m$-indescribable if for every $U\subseteq V_\kappa$  and
every $\Pi^n_m$ sentence $\phi$ (containing in prenex form $m$
alternations of $n$-th order quantifiers starting with $\forall$),
if $(V_\kappa,\in, U)\models\phi$, then there is $\alpha<\kappa$
such that $(V_\alpha,\in,\linebreak U\cap V_\alpha)\models\phi$.
The following is standard (see \cite[p. 297]{Je03} for a proof).

\begin{Thm} \label{T:Hanf}
{\rm (Hanf-Scott)} A cardinal $\kappa$ is weakly compact iff it is
$\Pi^1_1$-indescribable.
\end{Thm}

\begin{Def} \label{D:transfer}
{\rm (LZFC)} {\em A transitive model $M\models{\rm ZFC}$ is said
to be} $\Pi^1_1$-indescribable {\em if for every $U\in Def(M)$ and
every $\Pi^1_1$ sentence $\phi$, if $(M,\in, U, Def(M))\models
\phi$, then there is a transitive model $N\in M$ such that $U\cap
N\in Def(N)$ and $(N,\in, \linebreak U\cap N, Def(N))\models
\phi$.}
\end{Def}

In the above notation $Def(M)$, $Def(N)$ indicate the ranges for
the  second order quantifiers of $\phi$. $\Pi^1_1$-indescribability is first-order definable and absolute for transitive models. That is,  ``$M$ is
$\Pi^1_1$-indescribable'' iff $K\models\mbox{``$M$ is
$\Pi^1_1$-indescribable''}$ for any transitive model  $K$ such
that $M\in K$.

That $\Pi^1_1$-indescribable models (can be consistently assumed
to)  exist is a consequence of the following:

\begin{Prop} \label{P:rightgener}
{\rm (ZFC)} If $\kappa$ is weakly compact, then the model
$V_\kappa$ is $\Pi^1_1$-indescribable.
\end{Prop}

{\em Proof.} This is immediate from \ref{T:Hanf} and  lemma
\ref{L:flows} below. \telos

\begin{Lem} \label{L:flows}
{\rm (ZFC)} Let $\kappa$ be a $\Pi^1_1$-indescribable cardinal.
Then for every $U\in Def(V_\kappa)$,  and every $\Pi^1_1$ sentence
$\phi$ of ${\cal L}_2\cup\{\bf{S}\}$ (where ${\cal L}_2$ is ${\cal
L}$ augmented with second order variables and  $\bf{S}(\cdot)$ is
a unary predicate interpreted as $U$), if
$$(V_\kappa,\in,U,Def(V_\kappa))\models\phi,$$ then there
is $\alpha<\kappa$ such that $U\cap V_\alpha$ is (first-order)
definable in $(V_\alpha,\in)$ and $(V_\alpha,\in,U\cap
V_\alpha,Def(V_\alpha))\models\phi$.
\end{Lem}

{\em Proof.} Let $\kappa$ be  $\Pi^1_1$-indescribable. Let $U\in
Def(V_\kappa)$, and let $U=\{x\in
V_\kappa:V_\kappa\models\theta(x)\}$, for a first-order formula
$\theta$.  Let also $\phi=\forall X\psi(X)$ be a $\Pi^1_1$
sentence, where $\psi(X)$ has no second order variables. Suppose
$(V_\kappa,\in,U,Def(V_\kappa))\models\phi$. Set $\sigma=\forall
x({\bf S}(x)\leftrightarrow \theta(x))$. Then clearly
$(V_\kappa,\in,U,Def(V_\kappa))\models\sigma$ and $\sigma$ is
first-order. So
$$(V_\kappa,\in,U,Def(V_\kappa))\models \forall X\psi(X)
\wedge \sigma,$$ or equivalently
\begin{equation} \label{E:Def}
(V_\kappa,\in,U)\models (\forall X)(X\in Def(V_\kappa) \rightarrow
\psi(X))  \wedge \sigma.
\end{equation}
Now it is well-known that $Def(V_\kappa)$ is
$\Delta^1_1$-definable over $V_\kappa$.\footnote{Namely, $X\in
Def(V_\kappa):= (\exists \phi)(\forall x)(x\in X\leftrightarrow
Sat(\phi,x))$, where $Sat(\phi,x)$ is the $\Delta^1_1$
satisfaction predicate for first order formulas with parameters over $V_\kappa$.}
Therefore $(\forall X)(X\in Def(V_\kappa) \rightarrow \psi(X))
\wedge \sigma$ is $\Pi^1_1$ and hence, by $\Pi^1_1$-
indescribability of $\kappa$, there is $\alpha<\kappa$ such that
\begin{equation} \label{E:Def1}
(V_\alpha,\in,U\cap V_\alpha)\models (\forall X)(X\in
Def(V_\alpha) \rightarrow \psi(X)) \ \wedge \sigma.
\end{equation}
By the definition of $\sigma$,  $(V_\alpha,\in,U\cap
V_\alpha)\models \sigma$ implies that $U\cap V_\alpha=\{x\in
V_\alpha:V_\alpha\models\theta(x)\}$, that is, $U\cap V_\alpha\in
Def(V_\alpha)$. Further $(V_\alpha,\in,U\cap V_\alpha)\models
(\forall X)(X\in Def(V_\alpha) \rightarrow \psi(X))$ implies that
$$(V_\alpha,\in,U\cap V_\alpha, Def(V_\alpha))\models (\forall
X)\psi(X),$$ or $(V_\alpha,\in,U\cap
V_\alpha,Def(V_\alpha))\models\phi$, as required. \telos

\begin{Prop} \label{P:works}
If $M$ is a $\Pi^1_1$-indescribable model of {\rm ZFC} then $M$ is
$\alpha$-Mahlo for every $\alpha\in On^M$.
\end{Prop}

{\em Proof.} By induction on $\alpha$. Since $M$ is a model of
ZFC, it is $0$-Mahlo.  Let $\alpha=1$. We have to show that
$\{x\in M:(x,\in)\models{\rm ZFC}\}$ is stationary. Let $C\in
Club(M)$. There is a  first-order formula $\theta(x)$ such that
$x\in C\leftrightarrow M\models \theta(x)$. Let $\sigma=\forall
x({\bf S}(x)\leftrightarrow \theta(x))$. The fact that $C$ is a
club definable by $\theta(x)$ is expressed by writing
$$(M,\in,C,Def(M))\models \sigma \wedge ``\{x:{\bf S}(x)\} \
\mbox{is a club''}.$$  The sentence $\sigma \wedge ``\{x:{\bf
S}(x)\} \ \mbox{is a club''}$ is first-order so, by
$\Pi^1_1$-indescribability, there is $N\in M$, $N\models{\rm
ZFC}$, such that
$$(N,\in,C\cap N,Def(N))\models \sigma \wedge ``\{x:{\bf S}(x)\} \
\mbox{is a club''}.$$ This means that $\theta(x)$ defines  $C\cap
N$ in $N$ and   $C\cap N$ is a club of $N$. So if
$N_\alpha=V_\alpha^N$ for $\alpha\in N$,  we can pick by
induction, using Choice, sets $x_\alpha\in C\cap N$, $\alpha\in
M$, such that $N_\alpha \cup
(\bigcup_{\beta<\alpha}x_\beta)\subseteq x_\alpha$. If
$X=\{x_\alpha:\alpha\in N\}$, then clearly $X\in M$, $X\subseteq
C$ and $X$ is a chain. Therefore $\bigcup X\in C$. But $\bigcup
X=N$, so $N\in C$. It follows that the arbitrary club $C$ of $M$
contains a model $N\models{\rm ZFC}$.  Therefore $\{x\in
M:(x,\in)\models{\rm ZFC}\}$ is stationary in $M$.

Suppose $M$ is $(\alpha+1)$-Mahlo. Let $C\subseteq M$ be again a
club defined by $\theta(x)$ in $M$ and let $\sigma$ be as above.
Then
$$(M,\in,C,Def(M))\models \ \sigma \wedge \ \{x:{\bf S}(x)\} \
\mbox{is a club}\ \wedge \forall X(\mbox{$X$ is a
club}\rightarrow$$ $$\exists y(y\in X \wedge mahlo(\alpha,y)).$$
The last formula is $\Pi^1_1$ over $(M,\in,C,Def(M))$ and says
that $C$ is a club and that the definable set $\{x\in M:(x,\in) \
\mbox{is $\alpha$-Mahlo}\}$ is  a stationary set of $M$.  By
definition \ref{D:transfer}, there is $N\in M$, $N\models{\rm
ZFC}$, such that
$$(N,\in,C\cap N,Def(N))\models \ \sigma \wedge \ \{x:{\bf S}(x)\} \
\mbox{is a club}\ \wedge \forall X(\mbox{$X$ is a
club}\rightarrow$$ $$\exists y(y\in X \wedge mahlo(\alpha,y)).$$
This says that $C\cap N$ is a club of $N$ defined by $\theta(x)$
in $N$ and the set of  $\alpha$-Mahlo models contained in $N$ is a
stationary subset of $N$.  It follows that $N$ is
$(\alpha+1)$-Mahlo. Moreover, by the same argument as before, we
see that $N\in C$. So the arbitrary club $C$ of $M$ contains an
$(\alpha+1)$-Mahlo model. Therefore the set $\{x\in M: (x,\in) \
\mbox{is $(\alpha+1)$-Mahlo}\}$ is stationary in $M$, and hence
$M$ is $(\alpha+2)$-Mahlo.

Suppose $\alpha$ is limit and $M$ is $\alpha$-Mahlo. To show that $M$ is
$(\alpha+1)$-Mahlo the proof is essentially the same as before.

Finally, if $\alpha$ is limit and  the claim holds for all
$\beta<\alpha$, then, due to the definition of $\alpha$-Mahlo, the
claim holds for $\alpha$. \telos

\vskip 0.2in

{\sc Question.} What other large cardinal properties
(measurability, strong compactness, etc) can be adjusted to fit to
models of ZFC?

\section{Localizing extensions of ZFC}
In section 5 we have already considered extensions of $Loc({\rm
ZFC})$ of the form $Loc({\rm ZFC}+Loc({\rm ZFC}))$, $Loc({\rm
ZFC}+Loc({\rm ZFC}+Loc({\rm ZFC})))$, etc.   Here we shall
consider more general extensions, namely   localization principles
of the form $Loc({\rm ZFC}+\phi)$ for various sentences $\phi$
independent from ZFC. In order however for $Loc({\rm ZFC}+\phi)$
to make sense we must first assume that ${\rm ZFC}+\phi$ not only
is consistent but has a transitive model. So by analogy with the
axiom $TM({\rm ZFC})$ (``ZFC has a transitive model''), for every
such $\phi$ one has to accept
$$(TM({\rm ZFC}+\phi)) \hspace{.5\columnwidth minus .5\columnwidth}
 \exists x(Tr(x) \wedge (x,\in)\models{\rm
 ZFC}+\phi).
\hspace{.5\columnwidth minus .5\columnwidth} \llap{}$$ For several
natural sentences  like $V=L$, ${\rm CH}$,  $V\neq L$, $\neg {\rm CH}$ etc, it
is provable in ZFC (by usual  forcing techniques, constructible
sets, etc) that $TM({\rm ZFC})\rightarrow TM({\rm
ZFC}+\phi)$.\footnote{However one cannot prove in ZFC, if ZFC is consistent,  the implication $TM({\rm ZFC})\rightarrow TM({\rm ZFC}+TM({\rm
ZFC}))$, otherwise ${\rm ZFC}+TM({\rm ZFC})\vdash Con({\rm
ZFC}+TM({\rm ZFC}))$, contrary to G\"{o}del's  incompleteness. In
particular, ZFC does not prove, if it is consistent, that there is a  forcing extension $M[G]$
of a transitive model, that contains a transitive model of ZFC. } The same proof can be carried out (relativized) in
LZFC. Actually given a transitive model $M\models{\rm ZFC}$, there
is, by $Loc({\rm ZFC})$, a transitive $N$ such that $M\in N$. In
$N$ we can find a countable model $M'$ of ZFC and then extend it
by forcing, e.g. to a model $M'[G]$ of ${\rm ZFC}+\neg {\rm CH}$.

In the formulas occurring below as arguments in $Loc(\cdots)$ we allow the use of a constant ``$c$''. This is not a parameter, but ranges over  definable classes that are proved in ZFC to be sets (like  ${\cal P}(\omega)$, $\omega_1$, etc). Below we refer to such classes as ``terms''.  For the same reason ordinals occurring as parameters in formulas occurring as arguments in  $Loc(\cdots)$ are definable too.

Axioms $Loc({\rm ZFC}+\phi)$, though local in essence,  may have
global consequences for the universe $V$ itself. For example:

\begin{Lem} \label{P:allconstructible}
{\rm (LZFC)} Let $c$ be a term. Then
 $Loc( ZFC+V=L(c)) \rightarrow V=L(c)$.
\end{Lem}

{\em Proof.} Assume $Loc({\rm ZFC}+V=L(c))$ and  $V\neq L(c)$. Let
$a\in  V-L(c)$. Then  there is a transitive model $M$ of ZFC such
that $\{c,a\}\in M$ and $M\models  V=L(c)$. But then $a\in M=
 L(c)^M\subseteq L(c)$, a contradiction. \telos

\vskip 0.2in

More generally, given a set of sentences $\Gamma$, we may extend
LZFC to $${\rm LZFC}_\Gamma={\rm LZFC}+\{Loc({\rm
ZFC}+\phi):\phi\in \Gamma\} $$ and consider its  consistency  and
its consequences on  $V$. The following is a simple general fact
concerning  the consistency of ${\rm LZFC}_\Gamma$.

\begin{Prop} \label{P:compat}
If $\Gamma$ is a set of sentences such that
$\{\phi,\neg\phi\}\subseteq \Gamma$ for some $\Sigma_1^{\rm ZFC}$
or $\Pi_1^{\rm ZFC}$ sentence $\phi$, then ${\rm LZFC}_\Gamma$ is
inconsistent.
\end{Prop}

{\em Proof.} Let  $\phi$ be  a $\Sigma_1^{\rm ZFC}$ sentence (the
case of $\Pi_1^{\rm ZFC}$ is the same). This means that
$\phi\leftrightarrow \exists x\phi_1(x)$ holds in every model of
ZFC, for some $\Delta_0$ formula  $\phi_1$.   ${\rm ZFC}_\Gamma$
contains the axioms $Loc({\rm ZFC}+\phi)$ and $Loc({\rm
ZFC}+\neg\phi)$. By the first of them there is a transitive $M$
such that $M\models{\rm ZFC}+\phi$. Then $M\models\exists x
\phi_1(x)$, hence $M\models\phi_1(a)$ for some $a\in M$. By
$Loc({\rm ZFC}+\neg\phi)$ there is $N$ such that $a\in N$ and
$N\models{\rm ZFC}+\neg\phi$. Then $N\models \forall
x\neg\phi_1(x)$. But $N\models\phi_1(a)$, since $\phi_1$ is
absolute. A contradiction. \telos

\vskip 0.2in

Does any reasonable set $\Gamma$ affect the status of the axioms of Powerset, Separation,  Replacement, etc? (Remember that LZFC itself  is  compatible with ZFC).  The answer is positive for Powerset. We show  that if
$\Gamma$ contains the sentences ${\rm CH}$ and $\neg {\rm CH}$,  then ${\rm LZFC}_\Gamma$ refutes   Powerset.

Given a term $c$ and a transitive model $M$, let $c^M$ denote the relativization of  $c$ with respect to $M$. Let us call a term $c$ {\em stable} if for every transitive $M$, $c\subseteq M \Rightarrow c^M=c$. For instance ${\cal P}(\omega)$ and $H(\omega_1)$ are stable terms, while $\omega_1$ is not.

\begin{Prop} \label{P:contra}
 Let $c$ be a stable term. Then the theory
$${\rm LZFC}+Loc({\rm ZFC}+|c|=\omega_1)+Loc({\rm ZFC}+|c|\neq
\omega_1) +  \ \mbox{\rm ``$c$ exists''} +  ``{\cal P}(\omega) \ \mbox{\rm exists''} $$
is inconsistent.
\end{Prop}

{\em Proof.} Suppose the above mentioned theory  is consistent and let  $K$ be a model of it.  In $K$,  $c$ and ${\cal P}(\omega)$  are  sets. By $Loc({\rm ZFC}+|c|=\omega_1)$ and Pair, we can
pick a model $M\in K$ of ZFC such that $\{c,{\cal P}(\omega)\}\subset
M$ and $M\models |c|=\omega_1$. The last relation says that in $M$ there is a bijection $h:c^M\rightarrow \omega_1^M$. Since $c$ is stable, $c^M=c$. Also, since ${\cal P}(\omega)\in M$, by Lemma \ref{L:settle} (ii),  $\omega_1^M=\omega_1\in M$. Therefore  $c$ and $\omega_1$ are both absolute in $M$ and $M\models c\sim \omega_1$. Hence also $K\models
c\sim \omega_1$. Let $h:c\rightarrow \omega_1$ be a
bijection in $K$. By $Loc({\rm ZFC}+|c|\neq\omega_1)$ and
Pair, there is a model $N\in K$ such that
$\{c,h,{\cal P}(\omega)\}\subset N$ and $N\models |c|\neq
\omega_1$. Again, by stability $c^N=c$, and by \ref{L:settle} (ii), $\omega_1^N=\omega_1$. Hence $N\models c\not\sim \omega_1$. But this contradicts the fact that $N$ already contains a bijection $h:c\rightarrow \omega_1$.  \telos

\begin{Cor} \label{C:nopow}
(i) For every stable term $c$, the theory
$${\rm LZFC}+Loc({\rm ZFC}+|c|=\omega_1)+Loc({\rm ZFC}+|c|\neq
\omega_1) + \ \mbox{\rm ``$c$ exists''} + \ \mbox{\rm Powerset}$$
is inconsistent.

(ii) In particular, the theory
$${\rm LZFC}+Loc({\rm ZFC}+|{\cal P}(\omega)|=\omega_1)+Loc({\rm ZFC}+|{\cal P}(\omega)|\neq
\omega_1) + \mbox{\rm Powerset},$$
or, equivalently,
$${\rm LZFC}+Loc({\rm ZFC}+{\rm CH})+Loc({\rm ZFC}+\neg {\rm CH}) + \mbox{\rm Powerset}$$
is inconsistent.
\end{Cor}

{\em Proof.} (i) This follows immediately from \ref{P:contra},  if we replace ``${\cal P}(\omega)$ exists'' with the stronger Powerset.

(ii) In (i) above we set  $c={\cal P}(\omega)$, which is stable. Then Powerset implies ``${\cal P}(\omega)$ exists'' and the claim follows.

[It's worth  noting  that,  for this specific term $c={\cal P}(\omega)$, the claim can be alternatively proved (without appealing to \ref{P:contra})  as follows: Suppose ${\cal P}(\omega)$ is a set. Pick a model  $M$  such that ${\cal P}(\omega)\in M$ and $M\models |{\cal P}(\omega)|=\omega_1$. Then pick a model $N$ such that $M\in N$, hence $M\subseteq N$,  and $N\models |{\cal P}(\omega)|\neq \omega_1$. Then either $N\models |{\cal P}(\omega)|<\omega_1$, or $N\models |{\cal P}(\omega)|>\omega_1$. The first option is obviously false. So $N\models |{\cal P}(\omega)|>\omega_1$. But by Lemma \ref{L:Galois} (i),  $M\models |{\cal P}(\omega)|=\omega_1$ and $M\subseteq N$ imply $N\models |{\cal P}(\omega)|\leq \omega_1$, a contradiction.] \telos

\vskip 0.2in

In relation to clause (ii) of the last Corollary we  point out the following
(recall that NM is the assertion ``there is a natural model of ZFC'').

\begin{Prop} \label{P:conjunct}
The theory
$${\rm LZFC}+Loc({\rm ZFC}+{\rm CH})+Loc({\rm ZFC}+\neg {\rm CH})$$ is consistent relative to ${\rm ZFC}+{\rm NM}$.
\end{Prop}

{\em Proof.}  The proof is an easy strengthening of  that of Proposition \ref{P:PC} (iii). Let $M$ be a model of ${\rm ZFC}+{\rm NM}$.  Then $H^M(\omega_1)$ satisfies the theory in question. Indeed, if $M_\kappa=V_\kappa^M$ is a natural model of ZFC (in the sense of $M$), then by the proof of \ref{P:PC} (iii),    $H^M(\omega_1)\models Loc({\rm ZFC})$. So for every $x\in H^M(\omega_1)$,  $x$ belongs to  a countable transitive model  $N\in H^M(\omega_1)$. Now  every such model $N$ containing $x$ can be generically extended to countable transitive models $N_1,N_2$, satisfying CH and $\neg$CH, respectively. Since  $N_1,N_2$  also belong to $H^M(\omega_1)$,  $H^M(\omega_1)$ satisfies both $Loc({\rm ZFC}+{\rm CH})$ and $Loc({\rm ZFC}+\neg {\rm CH})$. \telos

\vskip 0.2in

Finally we have a variant of \ref{P:contra} from  which ``${\cal P}(\omega)$ exists'' has been dropped.

\begin{Prop} \label{P:Sepfrom}
Let $c$ be a stable term and let $\alpha\neq \beta$ be two distinct definable ordinals. Then  the theory
$${\rm LZFC}+Loc({\rm ZFC}+|c|=\omega_\alpha)+ Loc({\rm ZFC}+|c|=\omega_\beta) +\ \mbox{\rm ``$c$ exists''}$$
is inconsistent.
\end{Prop}

{\em Proof.} We work in the aforementioned theory and suppose $\alpha<\beta$. $c$ is a definable set, absolute for the models they contain it, hence by $Loc({\rm ZFC}+|c|=\omega_\alpha)$, there is a model $M$ of ZFC such that $c\in M$ and $M\models |c|=\omega_\alpha$. Then, by $Loc({\rm ZFC}+|c|=\omega_\beta)$, there is a model $N$ of ZFC such that $M\in N$ and $N\models |c|=\omega_\beta$. But  $M\subseteq N$ and  by Lemma \ref{L:Galois}   $M\models |c|=\omega_\alpha$ implies  $N\models |c|\leq \omega_\alpha$, i.e.,  $\beta\leq \alpha$, contrary to the   assumption  $\alpha<\beta$. \telos

\section{A digression: Standard compactness}
For reasons  explained in the introduction, one of the goals of this paper was to promote  transitive models of ZFC to the status of first class citizens of the universe of sets, especially by postulating their ``omnipresence''. In particular, whenever we talk about models in LZFC, we mean {\em transitive} models. Given that models in general is the stuff of the various notions of compactness, the confinement to transitive models induces  natural refinements of corresponding  compactness notions. Specifically, in ordinary compactness one infers the existence of a model for a set of sentences from the existence of models for its finite parts. A natural question arisen from this fact is the following: Can we infer the existence of a {\em transitive} model for a set  of sentences $\Sigma$ in a language ${\cal L}'$ extending the language ${\cal L}$ of set theory, from the existence of {\em transitive} models for certain parts of $\Sigma$? Although otherwise unrelated to the rest of the paper, this question is well-motivated by our insistence on transitive models and  shall be dealt with in this section. The question we formulated above prompts the following definition.

\begin{Def} \label{L:stcomp}
{\rm (ZFC)} {\em A cardinal $\kappa$ is said to be}  standard
compact {\em if for every set of sentences $\Sigma$ of a finitary
language  ${\cal L}'\supseteq {\cal L}$ such that
$|\Sigma|=\kappa$, if every set $A\subseteq \Sigma$ such that
$|A|<\kappa$ has a transitive model, then $\Sigma$ has a
transitive model.}
\end{Def}

A first negative result is that the standard version of the
classical compactness theorem is false.

\begin{Prop} \label{L:classic}
{\rm (ZFC or LZFC)} There is an ${\cal L}'\supseteq {\cal L}$ and  a countable  set  $\Sigma$  of  sentences of ${\cal L}'$ such that
every finite subset of $\Sigma$ has a transitive model, while
$\Sigma$ does not. Therefore $\omega$ is not standard compact.
\end{Prop}

{\em Proof.} Let ${\cal L}'={\cal L}\cup\{c_n:n\in\omega\}$, and
let $\Sigma=\{c_{n+1}\in c_n:n\in \omega\}$ be a set of sentences
of ${\cal L}'$. Then clearly every finite subset of  $\Sigma$ has
a transitive model, while $\Sigma$ does not. \telos

\vskip 0.2in

Next let us make sure that standard compact cardinals exist under
the assumption of mild large cardinals. Recall that one of the
equivalent definitions of a weakly compact cardinal is the
following: $\kappa$ is weakly compact if any set $\Sigma$ of
sentences of the infinitary language ${\cal L}_{\kappa,\kappa}$,
which uses at most $\kappa$ non-logical symbols and is
$\kappa$-satisfiable (i.e., every $A\subseteq \Sigma$ with
$|A|<\kappa$ is satisfiable), is satisfiable.

Recall also that by Mostowski's theorem \ref{L:Skolem} (ii), if
$E$ is a binary relation on $X$ such that (a) $E$ is well-founded
and (b) $(X,E)\models {\rm Ext}$, then there is a (unique)
transitive set $M$ such that $(X,E)\cong (M,\in)$. Ext is the
ordinary extensionality axiom, while well-foundedness is expressed
by a sentence of ${\cal L}_{\omega_1,\omega_1}$ as follows:
$${\rm Wf}:=\neg
(\exists_{n<\omega}x_n)(\bigwedge_{n<\omega}(x_{n+1}\in x_n)),$$
where $\exists_{n<\omega}x_n$ is an abbreviation of the infinite
block of quantifiers $\exists x_1\exists x_2\cdots\exists
x_n\cdots$. Every transitive set satisfies Ext and Wf. Conversely,
every ${\cal L}$-structure $(X,E)$ such that $(X,E)\models{\rm
Ext} \wedge {\rm Wf}$ is isomorphic to a transitive model. This is
a key fact by which we can prove the following:

\begin{Lem} \label{L:weaklyc}
Every weakly compact cardinal $\kappa>\omega$ is standard compact.
\end{Lem}

{\em Proof.} Let $\kappa>\omega$ be a weakly compact cardinal and
let $\Sigma$ be a set of sentences of ${\cal L}'\supseteq {\cal
L}$ such that $|\Sigma|=\kappa$. Suppose  that  every $A\subseteq
\Sigma$ with $|A|<\kappa$ has a transitive model. Let
$\Sigma'=\Sigma\cup\{{\rm Ext,Wf}\}$. $\Sigma'$ is  a set of
sentences of ${\cal L}_{\omega_1,\omega_1}$, and hence of ${\cal
L}_{\kappa,\kappa}$. From  the assumption about $\Sigma$ and the
remarks concerning Wf, every $A\subseteq \Sigma'$ with
$|A|<\kappa$ has a (transitive) model. By weak compactness of
$\kappa$, $\Sigma'$ has a  model $(X,E)$. Since this satisfies Ext
and Wf, it is isomorphic to a transitive model $M$. Thus
$\Sigma'$, and therefore $\Sigma$, has a transitive model. \telos

\begin{Prop} \label{P:defcard}
{\rm(ZFC)} $\omega_1$ is not standard compact. Similarly for
$\omega_n$, for every $n\in \omega$.
\end{Prop}

{\em Proof.}   Let ${\cal L}'={\cal L}\cup
\{c\}\cup\{\dot{\alpha}:\alpha\leq \omega_1\}$, where $c$ and
$\dot{\alpha}$ are constants. We shall find a $\Sigma$ that
refutes standard compactness of $\omega_1$. Let $\Sigma$ be the
set of the following sentences of ${\cal L}'$:

(1) $Ord(c)$, $Ord(\dot{\alpha})$, for all $\alpha\leq \omega_1$.

(2) $\dot{\alpha}< \dot{\beta}$, for all  $\alpha< \beta\leq
\omega_1$.

(3) $c>\dot{\alpha}$, for all $\alpha<\omega_1$.

(4) $c<\dot{\omega_1}$.

(5)  $\forall x(x<\dot{\omega_1}\rightarrow x \ \mbox{is
countable})$. \\
Clearly $|\Sigma|=\omega_1$. Let $A\subseteq \Sigma$ with
$|A|<\omega_1$. Pick some $V_\xi$ such that $\omega_1\in V_\xi$
and let $\dot{\alpha}^{V_\xi}=\alpha$ for all $\alpha\leq
\omega_1$. Then we easily see that  $V_\xi\models A$ for some
interpretation $c^{V_\xi}\in \omega_1$. On the other hand suppose
there is a transitive  structure $(K,\in)$ such that
$K\models\Sigma$. Although $K$ need not be a model of ZFC,
$K\models Ord(\dot{\alpha})$ clearly entails that $\dot{\alpha}^K$
is an ordinal.  In view of (2) the mapping
$\alpha\mapsto\dot{\alpha}^N$ is strictly increasing. Therefore
$\alpha\leq \dot{\alpha}^K$ for every $\alpha\leq \omega_1$, and
hence $\omega_1\leq \dot{\omega_1}^K$. By (5)   every $x\in
\dot{\omega_1}^K$ is countable, so in particular
$\dot{\omega_1}^K=\omega_1$. In view of this and  (3) and (4), we
have $\alpha\leq \dot{\alpha}^K<c^K<\dot{\omega_1}^K=\omega_1$,
hence $\alpha<c^K<\omega_1$, for all $\alpha<\omega_1$, which is a
contradiction.

In the case of $\omega_n$ we just need to replace ``$\alpha$ is
countable'' with the appropriate sentence defining $\omega_n$,
namely: ``$x$ is countable or  of cardinality next to countable,
or next to next to  countable or,...., or next$^n$ to
countable''. \telos

\vskip 0.2in

The property of weak compactness (as well as that of standard
compactness) contains  the  condition that  the cardinality  of
non-logical symbols (or the cardinality) of $\Sigma$  be $\leq
\kappa$. If we drop this condition we have the property of strong
compactness: $\kappa$ is strongly compact if  for every set
$\Sigma$ of sentences of ${\cal L}_{\kappa,\kappa}$, if $\Sigma$
is $\kappa$-satisfiable, then $\Sigma$ is satisfiable. An
equivalent definition (see \cite[p. 37] {Ka97}) is the following:

\begin{Def} \label{D:sc}
{\rm (ZFC) } {\em A cardinal $\kappa$ is} strongly compact {\em
if for any set $X$, every $\kappa$-complete filter on $X$ can be
extended to a $\kappa$-complete ultrafilter on $X$.}
\end{Def}

\begin{Prop} \label{P:nomizo}
{\rm (ZFC)} Let $\lambda>\omega$ be a  strongly compact cardinal.
Then every cardinal $\kappa\geq \lambda$ such that
$\kappa^{<\kappa}=\kappa$ is standard compact.
\end{Prop}

{\em Proof.} The proof is a variant of the proof of compactness by
the use of ultraproducts (see \cite[Cor. 4.1.11]{CK73}). Let
$\kappa\geq \lambda>\omega$, where $\lambda$ is strongly compact
and $\kappa^{<\kappa}=\kappa$. Let  $\Sigma$ be an infinite  set
of sentences of a language ${\cal L}'\supseteq {\cal L}$ such that
$|\Sigma|=\kappa$. Suppose each $A\subseteq \Sigma$ with
$|A|<\kappa$ has a transitive model. Since
$\kappa^{<\kappa}=\kappa$, there is an   enumeration
$\Sigma_\alpha$, $\alpha<\kappa$, of all subsets $A$ of $\Sigma$
with $|A|<|\Sigma|$. Pick and fix for each $\alpha<\kappa$ a
transitive model $M_\alpha\models\Sigma_\alpha$. For every
$\phi\in \Sigma$, let $\hat{\phi}=\{\beta<\kappa:\phi\in
\Sigma_\beta\}$. The family $E=\{\hat{\phi}:\phi\in \Sigma\}$ is
$\kappa$-complete, i.e., for every $\gamma<\kappa$ and every
$\{\hat{\phi_\beta}:\beta<\gamma\}\subseteq E$,
$\bigcap_{\beta<\gamma}\hat{\phi_\beta}\neq \emptyset$. This is
because for every $\gamma<\kappa$ and  every set
$\{\phi_\beta:\beta<\gamma\}$, there is a $\delta<\kappa$ such
that $\{\phi_\beta:\beta<\gamma\}=\Sigma_\delta$,  so $\delta\in
\bigcap_{\beta<\gamma} \hat{\phi_\beta}$. Thus the filter
$\bar{E}$ on $\kappa$ generated by $E$ is $\kappa$-complete. Also
$\bar{E}$ is free, otherwise  some $\alpha$ would be in all
$\hat{\phi}$, $\phi\in \Sigma$, hence $\Sigma=\Sigma_\alpha$,
which is false.  One can see as in \cite[Cor. 4.1.11]{CK73} that
if $D$ is any ultrafilter on $\kappa$ extending $\bar{E}$, then
$\Pi_{\alpha<\kappa} M_\alpha/D\models\Sigma$. Namely for every
$\phi\in \Sigma$,
$$\{\alpha<\kappa:M_\alpha\models\phi\}\supseteq
\{\alpha<\kappa:\phi\in \Sigma_\alpha\}=\hat{\phi}\in D,$$ so
$\Pi_{\alpha<\kappa} M_\alpha/D\models\phi$ by the fundamental
theorem of ultraproducts. It suffices to choose the ultrafilter
$D\supseteq \bar{E}$ so that $\Pi_{\alpha<\kappa} M_\alpha/D$ be
(isomorphic to) a transitive model.  Now $\bar{E}$ is a  $\kappa$-complete
filter and hence $\lambda$-complete since $\lambda\leq \kappa$. But $\lambda$ is strongly compact, so $\bar{E}$ can be extended to a $\lambda$-complete ultrafilter $D$.  Since $\lambda>\omega$ and
every $M_\alpha$ is transitive,  the ultraproduct
$\Pi_{\alpha<\kappa} M_\alpha/D$ is well-founded. Therefore
$\Pi_{\alpha<\kappa} M_\alpha/D$ is isomorphic to a transitive
$(N,\in)$. Then $(N,\in)\models\Sigma$ as required. \telos

\vskip 0.2in

It follows from the last result that, unless strongly compact cardinals are inconsistent, it is consistent to have standard compact cardinals which are accessible, singular and even successor cardinals.

\vskip 0.2in

{\bf Acknowledgement.} Many thanks to the anonymous referee for carefully  checking the manuscript, pointing out some serious flaws and  suggesting a lot of  other  improvements.

\end{document}